# RESAMPLING METHODS FOR SPATIAL REGRESSION MODELS UNDER A CLASS OF STOCHASTIC DESIGNS[1]

By S. N. Lahiri and Jun Zhu

*Iowa State University and University of Wisconsin–Madison*

In this paper we consider the problem of bootstrapping a class of spatial regression models when the sampling sites are generated by a (possibly nonuniform) stochastic design and are irregularly spaced. It is shown that the natural extension of the existing block bootstrap methods for grid spatial data does not work for irregularly spaced spatial data under nonuniform stochastic designs. A variant of the blocking mechanism is proposed. It is shown that the proposed block bootstrap method provides a valid approximation to the distribution of a class of $M$-estimators of the spatial regression parameters. Finite sample properties of the method are investigated through a moderately large simulation study and a real data example is given to illustrate the methodology.

**1. Introduction.** Irregularly spaced spatial data occur frequently in many applications and seem to be the rule rather than the exception (cf. [15]). Distributional properties of estimators and test statistics based on irregularly spaced spatial data are typically influenced by intricate interactions among the design density, the spatial asymptotic structure and the spatial correlation structure of the underlying process, and have complicated forms. Resampling methods can play a very useful role in the statistical analysis of such spatial data because these methods do not require analytical derivations of the exact forms of such interactions and their effects on the (asymptotic) distribution of a statistic. Although work on resampling methods for spatial data (cf. [13]) began at about the same time as that for temporal processes (cf. [7]), much less seems to be known about resampling methods for spatial data that are irregularly spaced. In a pioneering work, Politis, Paparoditis and Romano [29] developed a subsampling method for irregularly spaced

Received December 2003; revised March 2005.
[1]Supported in part by NSF Grants DMS-00-72571 and DMS-03-06574.

spatial data generated by a *homogeneous* Poisson process. Politis, Paparoditis and Romano [30] formulated a version of the spatial block bootstrap, also under the same framework. In this paper we consider the problem of formulating a unified spatial block bootstrap method for irregularly spaced spatial data in a more general framework where the data-sites may have a *nonuniform* concentration across the sampling region and where the data generating mechanism may have a nontrivial *infill sampling* component to it. Both these features are important in applications involving irregularly spaced spatial data, as (i) the total number of sampling sites and the volume of the sampling region may not be of a comparable order of magnitude (as required by a homogeneous Poisson process) and (ii) a higher concentration of sampling sites may be located around the "hot-spots" (leading to a nonuniform density of sampling sites).

To highlight the peculiarities associated with formulation of a valid blocking mechanism under a nonuniform sampling design, we consider the following example. Let $R$ be a sampling region in the plane. First suppose that the sampling sites are generated by a square-grid, for example, the integer grid $\mathbb{Z}^2$, where $\mathbb{Z} = \{0, \pm 1, \ldots\}$ denotes the set of all integers. For this case, Politis and Romano [31] give a formulation of the block bootstrap, where (overlapping) blocks are formed by considering translations of a (square-)template by points on the grid. Thus, in the regular grid case, one uses the data-locations themselves to define the blocks. A random sample of these blocks then yields the bootstrap observations.

Now consider the case where the data-sites are generated by a stochastic design and hence are irregularly spaced. In analogy to the regular grid case, one may form the blocks by using translations of a given (square-)template of a suitable size by the sampling sites themselves [see Figure 1(b)] and resample from these blocks to generate the "bootstrap observations." We call this version of the spatial block bootstrap the *data-site-shifted block bootstrap* method or the DSSBB method, in short. Although the DSSBB method is a "natural" extension of the standard spatial block bootstrap method for regularly spaced grid data to the irregularly spaced data case, it turns out that the DSSBB method may fail when the spatial sampling density is nonuniform. In Section 4 we construct an example of a nonuniform spatial sampling density and show that the DSSBB estimator of the variance of the sample mean is inconsistent for a wide range of block sizes. As a result, the DSSBB is not suitable for irregularly spaced spatial data when the design density is nonuniform.

As a remedy, we also present an alternative formulation of the blocking mechanism and establish its validity for irregularly spaced spatial data, allowing nonuniformity of the spatial sampling design and also allowing nonstationarity of the underlying spatial process. The alternative formulation introduces an extraneous regular grid and defines the blocks by considering translates of a suitable template by points on this grid [see Figure 1(c)]. Although the alternative blocking mechanism is less natural compared to the version based on the data-site-shifted blocks, it circumvents the problems associated with the natural blocking mechanism of the DSSBB method. We refer to this new version of the block bootstrap method as the *grid-based block bootstrap* method or the GBBB method, in short.

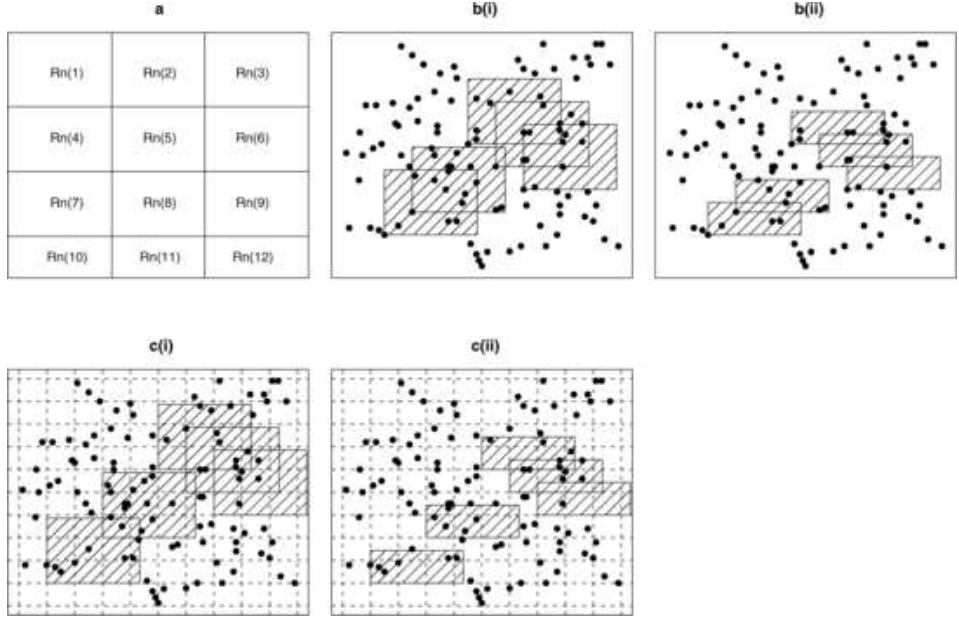

Fig. 1. (a) *Sampling region and the partition* $R_n = \bigcup_{\mathbf{k} \in \mathbf{K}_n} R_n(\mathbf{k})$ [*cf.* (4.3)]. *For simplicity, we write* $R_n(\mathbf{k})$ *as* $R_n(1), \ldots, R_n(12)$ *in the figure.* [b(i)] *A few representative translates of a complete block, say* $R_n(1)$, *with the lower left points on the sampling sites. These translates are used to construct the DSSBB version of the spatial process on the complete blocks.* [b(ii)] *A few representative translates of an incomplete block, say* $R_n(10)$, *with the lower left points on the sampling sites. These translates are used to construct the DSSBB version of the spatial process on the incomplete blocks.* [c(i)] *A few representative translates of a complete block, say* $R_n(1)$, *with the lower left points on the introduced regular grid. These translates are used to construct the GBBB version of the spatial process on the complete blocks.* [c(ii)] *A few representative translates of an incomplete block, say* $R_n(10)$, *with the lower left points on the introduced regular grid. These translates are used to construct the GBBB version of the spatial process on the incomplete blocks.*

In this paper we investigate properties of the GBBB method for a class of spatial regression models under a stochastic design that is driven by a collection of independent and identically distributed (i.i.d.) random vectors with a possibly nonuniform density. The spatial asymptotic structure adopted here allows "infilling" of any given subregion of the sampling design and thus, provides a more flexible framework than the homogeneous Poisson process formulation that is commonly assumed for modeling irregularly spaced spatial data. The main theoretical results of the paper show that under mild regularity conditions, the GBBB method provides a consistent estimator of the asymptotic variance of a class of $M$-estimators of the regression parameters. We also show that the GBBB provides a valid approximation to the distributions of these $M$-estimators. Thus, the GBBB provides a unified bootstrap method for irregularly spaced spatial data that gives a valid approximation, irrespective of nonuniformity of the data-sites, even in presence of infill sampling. Finite sample properties of the GBBB method are studied in a moderately large simulation study. As an illustration of the proposed methodology, we also consider a real data example involving the grazing patterns of elks and landscape characteristics in northern Wisconsin, where elks are of high conservation interest.

Block bootstrap methods for time-series data and for grid spatial data have been put forward by Hall [13], Carlstein [7], Künsch [17], Liu and Singh [26], Politis and Romano [31, 32], Zhu and Morgan [41], Zhu and Lahiri [39], among others. Subsampling methods for grid spatial data have been developed by Possolo [35], Politis and Romano [33], Sherman and Carlstein [37], Sherman [36], Garcia-Soidan and Hall [11], Lahiri [20], Lahiri, Kaiser, Cressie and Hsu [23] and Zhu, Lahiri and Cressie [40]. Optimal block sizes for the spatial resampling method have been given by Nordman and Lahiri [27]. Lee and Lahiri [25] employ spatial resampling to develop a least squares method for estimation of covariance parameters. Bootstrap and subsampling methods for irregularly spaced spatial data have been formulated by Politis, Paparoditis and Romano [29, 30]. See Chapter 12 of [21] for more on resampling methodology for spatial data.

The rest of the paper is organized as follows. In Section 2 we introduce the stochastic design and the spatial asymptotic framework. In Section 3 we describe the spatial regression model and establish limit distribution theory for a class of $M$-estimators under the stochastic design of Section 2. These results appear to be new and may be of some independent interest. Among other things, the results of Section 3 highlight the effects of the nonuniform sampling design and of the spatial asymptotic structure on the large sample distribution of the $M$-estimators of the regression parameter. In Section 4 we present the DSSBB method and give an example to show its inadequacy under nonuniform stochastic designs. In Section 5 we describe the GBBB method and establish theoretical results on the validity of the GBBB approximation for the distribution of the $M$-estimators. Results from a moderately large simulation study and the real data example are presented in Section 6 and Section 7, respectively. For clarity of exposition, proofs of the main theoretical results are presented in Sections 8, 9 and 10.

**2. The spatial sampling design.**

2.1. *Sampling regions.* Let $U_0$ be an open connected subset of $(-1/2, 1/2]^d$ containing the origin and let $R_0$ be a Borel set satisfying $U_0 \subset R_0 \subset \bar{U}_0$, where for any set $A \subset \mathbb{R}^d$, $\bar{A}$ denotes its closure. We regard $R_0$ as a "prototype" of the sampling region $R_n$. Let $\{\lambda_n\}$ be a sequence of positive numbers such that $n^\varepsilon/\lambda_n \to 0$ as $n \to \infty$ for some $\varepsilon > 0$. We assume that the sampling region $R_n$ is obtained by "inflating" the set $R_0$ by the scaling factor $\lambda_n$ (cf. [23]), that is,

$$R_n = \lambda_n R_0. \tag{2.1}$$

Since the origin is assumed to lie in $R_0$, the shape of $R_n$ remains the same for different values of $n$. To avoid pathological cases, we assume that for any sequence of real numbers $\{a_n\}$ with $a_n \to 0+$ as $n \to \infty$, the number of cubes of the lattice $a_n \mathbb{Z}^d$ that intersect both $R_0$ and $R_0^c$ is $O((a_n)^{-(d-1)})$ as $n \to \infty$. The boundary condition on $R_0$ holds for most regions $R_n$ of practical interest, including common convex subsets of $\mathbb{R}^d$, such as spheres, ellipsoids, polyhedrons, as well as for many nonconvex star-shaped sets in $\mathbb{R}^d$. (Recall that a set $A \subset \mathbb{R}^d$ is called star-shaped if for any $x \in A$, the line segment joining $x$ to the origin lies in $A$.) The latter class of sets may have fairly irregular shape. See, for example, [37] and [20] for more details.

2.2. *The stochastic sampling design.* Let $f(x)$ be a probability density function on $R_0$ and let $\{\mathbf{X}_n\}$ be a sequence of independent and identically distributed (i.i.d.) random vectors with probability density function $f(x)$ such that $\{\mathbf{X}_n\}$ and $\{Z(\mathbf{s}): \mathbf{s} \in \mathbb{R}^d\}$ are defined on a common probability space $(\Omega, \mathcal{F}, P)$ and are independent. We assume that the sampling sites $\mathbf{s}_1, \ldots, \mathbf{s}_n$ are obtained from a realization $\mathbf{x}_1, \ldots, \mathbf{x}_n$ of the random vectors $\mathbf{X}_1, \ldots, \mathbf{X}_n$, by the relation

$$\mathbf{s}_i = \lambda_n \mathbf{x}_i, \qquad 1 \le i \le n. \tag{2.2}$$

Since $\mathbf{x}_1, \ldots, \mathbf{x}_n$ take values in $R_0$, $\mathbf{s}_1, \ldots, \mathbf{s}_n$ take values in the entire sampling region $R_n \equiv \lambda_n R_0$. A similar spatial stochastic design has been used by Hall and Patil [15] in the context of nonparametric estimation of the autocovariance function of a spatial process. Note that by the strong law of large numbers (SLLN), the number of sampling sites lying in any given subregion $A$ of $R_n$ is asymptotically equivalent to $nP(\lambda_n X_1 \in A) = n \int_{\lambda_n^{-1} A} f(x)\,dx$. As a result, for a density function $f(x)$ that is not a constant function on $R_0$, one may have different expected numbers of sampling sites in two distinct subregions $A_1$ and $A_2$ of $R_n$ having the same volume. Thus, this formulation allows us to model irregularly spaced sampling sites that possibly have *nonuniform* concentration across different parts of the sampling region $R_n$. In contrast, given the sample size $n$, the standard approach of modeling irregularly spaced sampling sites using a homogeneous Poisson point process allows the sampling sites to have *only* the uniform distribution over the sampling region. A second important feature of the stochastic sampling design is that it allows the sample size $n$ and the volume of the sampling region $R_n$ to grow at different rates. For a positive design density $f$, when the sample size $n$ grows at a rate faster than $\text{vol}(R_n)$, the volume of $R_n$, the ratio of the expected number of sampling sites in a given subregion $A$ of $R_n$ to $\text{vol}(A)$ tends to infinity. This corresponds to the case of *infill* sampling in the stochastic design case (cf. [8], Chapter 5, [22]) and it is of interest in geostatistical and environmental monitoring applications (cf. [8, 23]). Thus, both these features of the present stochastic spatial sampling design offer some additional flexibility in handling irregularly spaced sampling sites over the standard homogeneous Poisson process formulation and are important for practical applications.

In the next section we describe a class of spatial regression models and a class of $M$-estimators of the regression parameter vector under the spatial sampling design described above.

### 3. M-estimation in spatial regression models.

3.1. *A class of spatial regression models.* We consider the spatial regression model

$$(3.1) \qquad Y(\mathbf{s}) = \mathbf{w}(\mathbf{s})'\beta + Z(\mathbf{s}), \qquad \mathbf{s} \in \mathbb{R}^d,$$

where $\mathbf{w}(\mathbf{s}): \mathbb{R}^d \to \mathbb{R}^p$ is a nonrandom weight function, $\beta \in \mathbb{R}^p$ is the vector of regression parameters, $\{Z(\mathbf{s}) : \mathbf{s} \in \mathbb{R}^d\}$ is a zero-mean stationary random field (r.f.) on $\mathbb{R}^d$, $p, d \in \mathbb{N} \equiv \{1, 2, \ldots\}$ and where we denote the transpose of a matrix $A$ by $A'$. The weight function $\mathbf{w}(\mathbf{s})$ in (3.1) may depend on a set of covariates $x_1(\mathbf{s}), \ldots, x_q(\mathbf{s})$ as well as on the spatial location $\mathbf{s}$. Thus, the general form of $\mathbf{w}(\mathbf{s})$ is

$$\mathbf{w}(\mathbf{s}) = \gamma(x_1(\mathbf{s}), \ldots, x_q(\mathbf{s}), \mathbf{s})$$

for some function $\gamma : \mathbb{R}^{q+d} \to \mathbb{R}^p$. Next suppose that the $Y$-process is observed at finitely many locations $\mathbf{s}_1, \ldots, \mathbf{s}_n$ in the sampling region $R_n \subset \mathbb{R}^d$, generated by the stochastic sampling design of Section 2, and let $\psi : \mathbb{R} \to \mathbb{R}$ be a Borel-measurable function satisfying

$$(3.2) \qquad E\psi(Z(\mathbf{0})) = 0.$$

We define an $M$-estimator $\hat{\beta}_n$ of $\beta$ based on $\{(Y(\mathbf{s}_i), \mathbf{w}(\mathbf{s}_i)')' : i = 1, \ldots, n\}$ as a measurable solution to the equation (in $\mathbf{t} \in \mathbb{R}^p$)

$$(3.3) \qquad \sum_{i=1}^n \mathbf{w}(\mathbf{s}_i) \psi(Y(\mathbf{s}_i) - \mathbf{w}(\mathbf{s}_i)'\mathbf{t}) = 0.$$

Note that $\psi_0(x) = x$, $x \in \mathbb{R}$, gives the least squares estimator of $\beta$.

Distributional properties of $M$-estimators $\hat{\beta}_n$ depend on the population characteristics of the underlying spatial process as well as on the spatial sampling design in a nontrivial manner. Here we establish large sample distributional properties of the $M$-estimator $\hat{\beta}_n$ under the stochastic design of Section 2.

3.2. *Asymptotic distribution of M-estimators.* To state the results, we need to introduce some notation. For $\mathbf{x} = (x_1, \ldots, x_d)' \in \mathbb{R}^d$, write $\|\mathbf{x}\|_1 \equiv |x_1| + \cdots + |x_d|$ and $\|\mathbf{x}\| = (x_1^2 + \cdots + x_d^2)^{1/2}$, respectively, for the $\ell^1$ and $\ell^2$ norms on $\mathbb{R}^d$. For any set $A \subset \mathbb{R}^d$, let $\text{vol}(A)$ denote the volume (i.e., the Lebesgue measure) of $A$ and let $|A|$ denote the size (i.e., the number of elements) of $A$. For $y \in \mathbb{R}$, write $y_+ = \max\{y, 0\}$. Let $\mathbb{1}(\cdot)$ denote the indicator function and let $C$, $C(\cdot)$ denote generic (nonrandom) positive constants, depending on their arguments (if any), but not on $n$. Unless explicitly stated, limits in order symbols are taken by letting $n \to \infty$.

Next we define the strong-mixing coefficient for the r.f. $\{Z(\cdot)\}$. Let $\mathcal{F}_Z(T) = \sigma\langle Z(\mathbf{s}) : \mathbf{s} \in T\rangle$ be the $\sigma$-field generated by the variables $\{Z(\mathbf{s}) : \mathbf{s} \in T\}$, $T \subset \mathbb{R}^d$. For any two subsets $T_1$ and $T_2$ of $\mathbb{R}^d$, let $\tilde{\alpha}(T_1, T_2) = \sup\{|P(A \cap B) - P(A)P(B)| : A \in \mathcal{F}_Z(T_1),\ B \in \mathcal{F}_Z(T_2)\}$, where $d(T_1, T_2) = \inf\{\|\mathbf{x} - \mathbf{s}\|_1 : \mathbf{x} \in T_1,\ \mathbf{s} \in T_2\}$. Then the strong-mixing coefficient of the r.f. $\{Z(\cdot)\}$ is defined as

$$(3.4) \quad \alpha(a; b) = \sup\{\tilde{\alpha}(T_1, T_2) : d(T_1, T_2) \geq a,\ T_1, T_2 \in \mathcal{R}(b)\}, \quad a > 0,\ b > 0,$$

where $\mathcal{R}(b)$ is the collection of all finite disjoint unions of cubes in $\mathbb{R}^d$ with a total volume not exceeding $b$. Note that the supremum in the definition of $\alpha(a; b)$ is taken over sets $T_1, T_2$ that are *bounded*. This restriction is important for $d \geq 2$, as shown by Bradley [5, 6]. For simplifying the exposition, we further assume [cf. [9]] that there exist constants $C, \tau_1 \in (0, \infty)$ and $\tau_2 \in [0, \infty)$ such that

$$(3.5) \quad \alpha(a; b) \leq C \cdot a^{-\tau_1} b^{\tau_2} \quad \text{for all} \quad a \geq 1,\ b \geq 1,$$

and that $\tau_2 = 0$ for $d = 1$. Thus, for $d \geq 2$, the strength of dependence increases with the volumes of the sets $T_1$ and $T_2$. We shall specify the exact conditions on $\tau_1, \tau_2$ in the statements below.

CONDITIONS.

(C.1) There exists a sequence of nonsingular matrices $\{\Lambda_n\}$ such that

$$(3.6) \quad \Lambda_n^{-1}[E\mathbf{w}(\lambda_n \mathbf{X}_1)\mathbf{w}(\lambda_n \mathbf{X}_1)'](\Lambda_n^{-1})' \to H \quad \text{as } n \to \infty,$$

and for all $\mathbf{h} \in \mathbb{R}^d$,

$$(3.7) \quad \Lambda_n^{-1}\left[\int \mathbf{w}(\lambda_n\mathbf{x}+\mathbf{h})\mathbf{w}(\lambda_n\mathbf{x})'f^2(\mathbf{x})\,d\mathbf{x}\right](\Lambda_n^{-1})' \to Q(\mathbf{h}) \qquad \text{as } n \to \infty,$$

where $H$ is a positive definite matrix and where $Q(\cdot)$ is a $p \times p$ matrix-valued function on $\mathbb{R}^d$.

(C.2) $\psi$ is differentiable and its first derivative $\psi'$ satisfies a Lipschitz condition of order $\gamma \in (0,1]$.

(C.3) There exist a $\delta > 0$ and $r \in \mathbb{N}$ such that:

(i) $E|\psi(Z(\mathbf{0}))|^{2r+\delta} < \infty$, $E|\psi'(Z(\mathbf{0}))|^{2r+\delta} < \infty$.
(ii) $\tau_1 > (2r-1)(2r+\delta)d/\delta$, where $\tau_1$ is as in (3.5).
(iii) $0 \leq \tau_2 < \frac{\tau_1-d}{4d}$, if $d \geq 2$.
(iv) $\int Q(\mathbf{h})\sigma_\psi(\mathbf{x})\,d\mathbf{x}$ is positive definite, where $\sigma_\psi(\mathbf{x}) = E\psi(Z(\mathbf{0}))\psi(Z(\mathbf{x}))$.
(v) $\chi_0 = E\psi'(Z(\mathbf{0})) \neq 0$.

(C.4) For some $a_0 \in [0, 1/8)$ and for some $\varepsilon_n \downarrow 0$,

$$m_{0n}^2 \equiv \sup\{\|\Lambda_n^{-1}\mathbf{w}(\mathbf{s})\| : \mathbf{s} \in R_n\} = O(\min\{n^{a_0}, \varepsilon_n(\log n)^{-2}\lambda_n^{(\tau_1-d)/4\tau_1}\}).$$

(C.5) The probability density $f(\cdot)$ of $\mathbf{X}_1$ is positive and continuous on $\bar{R}_0$.

We briefly comment on the conditions. The first part of condition (C.1) [cf. (3.6)] is formulated to add flexibility in the definition of the normalizing matrix for the $M$-estimator $\beta$. One possible choice of $\Lambda_n$ is given by $[E\mathbf{w}(\lambda_n\mathbf{x}_1)\mathbf{w}(\lambda_n\mathbf{x}_1)']^{1/2}$, in which case we may take $H = I_p$, the identity matrix of order $p$. However, this choice may not be the most convenient for the verification of (3.7), which is a version of the well-known Grenander condition in our spatial regression context. See [2, 12, 22] for more details about condition (3.7). Condition (C.2) is a set of smoothness conditions on the score function $\psi$ and is satisfied by many score functions $\psi$, such as the function $\psi_0(x) = x$, $x \in \mathbb{R}$, giving the least squares estimator of $\beta$, and by the log-pseudolikelihood function of certain Markov r.f. error processes. For certain *nondecreasing* score functions, the limit distributions of $M$-estimators $\hat{\beta}_n$ are derived in [24] using a different method of proof. Conditions (C.3)–(C.5) are used for deriving the limit distribution of the $M$-estimator $\hat{\beta}_n$ under the given spatial sampling design and are similar to the condition imposed in [22]. The choice of $r$ in (C.3) will be specified in the statements of the results below. For example, for establishing the asymptotic normality of the $M$-estimator (cf. Theorem 1), we used condition (C.3) with $r = 1$. Note that corresponding to a bounded $\psi$, the requirement on the strong-mixing coefficient in Theorem 1 reduces to $\alpha(a;b) \leq Ca^{-(d+\varepsilon)}b^{1/\varepsilon}$, $a > 1$, $b > 1$ for $d \geq 2$, and to $\alpha(a;\infty) \leq Ca^{-1-\varepsilon}$, $a > 1$ for $d = 1$, for some arbitrarily small $\varepsilon \in (0,1]$. The latter is close to the optimal rate for $d = 1$.

Next define the possible asymptotic covariance matrices of the normalized $\hat{\beta}_n$ by

(3.8) $$\Xi_c = \chi_0^{-2} H^{-1} \Sigma_c H^{-1} \qquad \text{for } c \in (0, \infty],$$

where

(3.9) $$\Sigma_\infty = \int \sigma_\psi(\mathbf{x}) Q(\mathbf{x}) \, d\mathbf{x} \quad \text{and} \quad \Sigma_c = c^{-1} H \sigma_\psi(\mathbf{0}) + \Sigma_\infty$$
$$\text{for } c \in (0, \infty).$$

Also, let $\mathcal{L}(\mathbf{W}|\mathcal{X})$ denote the conditional distribution of a random vector $\mathbf{W}$ given $\mathcal{X} \equiv \sigma \langle \{\mathbf{X}_n : n \geq 1\} \rangle$, and let $\varrho$ denote the Prohorov metric (cf. [4, 28]) on the set of all probability measures on $\mathbb{R}^p$. Then we have the following result.

THEOREM 1. *Suppose that conditions* (C.1)–(C.5) *hold with* $r = 1$. *Also, suppose that* $n/\lambda_n^d \to c$ *for some* $c \in (0, \infty]$. *Then*

$$\lambda_n^{d/2} \Lambda_n'(\hat{\beta}_n - \beta) \xrightarrow{d} N(0, \Xi_c) \qquad \text{for almost all } \mathbf{X}_1, \mathbf{X}_2, \ldots, (P_{\mathbf{X}}),$$

*that is,*

(3.10) $\varrho(\mathcal{L}(\lambda_n^{d/2} \Lambda_n'(\hat{\beta}_n - \beta)|\mathcal{X}), N(\mathbf{0}, \Xi_c)) \to 0 \qquad \text{as } n \to \infty, \ a.s. \ (P_{\mathbf{X}}),$

*where* $\Xi_c, c \in (0, \infty]$, *are as defined in* (3.8) *and* (3.10).

Theorem 1 shows that the $M$-estimator $\hat{\beta}_n$ is asymptotically normal for almost all realizations of the sampling design random vectors $\mathbf{X}_1, \mathbf{X}_2, \ldots$, with its asymptotic covariance matrix depending on the relative growth rates of the sample size $n$ and the volume, $|R_0|\lambda_n^d$, of the sampling region $R_n$. Let $c \equiv \lim_{n \to \infty} n/\lambda_n^d$. When $c \in (0, \infty)$, the sample size $n$ and the volume of the sampling region $R_n$ grow at the same rate, and this corresponds to the "pure increasing domain asymptotic structure" (cf. [8, 22]) in the present setup. On the other hand, when $c = \infty$, the sample size grows at a faster rate than the volume of $R_n$ and therefore, any given subregion of $R_n$ of unit volume may contain an unbounded number of sampling sites as $n \to \infty$. Thus, for $c = \infty$, the sampling region $R_n$ is subjected to infill sampling and we get a "mixed increasing domain asymptotic structure" with a nontrivial infill component. Theorem 1 establishes asymptotic normality of the $M$-estimators $\hat{\beta}_n$ under both types of spatial asymptotic structures. Furthermore, in view of (3.8)–(3.10), Theorem 1 also shows that the infill component leads to a reduction in the asymptotic variance of $\hat{\beta}_n$ in the mixed case, as the positive definite matrix $c^{-1} H \sigma(\mathbf{0})$ drops out for the "$c = \infty$" case.

REMARK 3.1. From the proof of Theorem 1, it follows that under the regularity conditions of Theorem 1, the estimating equation (3.3) admits a sequence of solutions $\{\hat{\beta}_n\}$ such that for some constant $C \in (0, \infty)$,

(3.11) $\quad P_{\cdot|\mathbf{X}}(\|\hat{\beta}_n - \beta\| > C\lambda_n^{-d/2}(\log n)^2) = o(1) \qquad \text{as } n \to \infty, \text{ a.s. } (P_\mathbf{X}),$

where $P_{\cdot|\mathbf{X}}(\cdot)$ denotes the conditional probability given $\mathcal{X}$. In particular, by the bounded convergence theorem, this sequence of solutions is "consistent" in the sense that it converges to $\beta$ in (the unconditional) probability. When the solution to (3.3) is unique, the unique solution also satisfies (3.11) and therefore is consistent for $\beta$. However, if (3.3) has multiple solutions, Theorem 1 applies to any sequence $\{\hat{\beta}_n\}$ that satisfies (3.11). In practice, when (3.3) has multiple solutions, one may identify such a sequence by considering the solution that is closest to a $\lambda_n^{d/2}$-consistent estimator (such as the least squares estimator) of $\beta$.

Relation (3.10) shows that the limit distributions of the $M$-estimators depend in a very complicated way on the correlation structure of the error process, on the regression function and on the spatial sampling density $f$. In view of the complicated form of the asymptotic variance and its dependence on the value of the infill parameter $c$, which is typically unknown in practice, it is important to develop resampling methods that automatically adjust themselves to the different forms of the sampling distribution of normalized $\hat{\beta}_n$ under various combinations of these factors and produce valid approximations in all cases.

In the next section we describe the natural extension of the standard block bootstrap method, the DSSBB, and give an example to show that it may fail if the design density is nonuniform. In Section 5 we describe the GBBB method and establish its validity for approximating the sampling distribution of $M$-estimators in the *full generality* of the spatial sampling framework of Theorem 1, including nonuniform densities $f$ and both types of spatial asymptotic structures.

**4. The DSSBB method and its properties.** For simplicity of exposition, we describe the DSSBB method only for $d$-dimensional *cubic* sampling regions for a *stationary* r.f. A description of the method for a sampling region of a general shape and for a spatial process satisfying (3.1) can be given by routinely modifying the description of the GBBB method given in the next section. Hence, for the rest of this section, suppose that $R_0 = (-1/2, 1/2)^d$ is the prototype set for the sampling region

(4.1) $\qquad\qquad R_n = \lambda_n R_0 = \left(-\frac{\lambda_n}{2}, \frac{\lambda_n}{2}\right)^d,$

and that $\{Y(\mathbf{s}) : \mathbf{s} \in \mathbb{R}^d\}$ is a stationary r.f. [i.e., $\beta = \mathbf{0}$ in (3.1)].

4.1. *The DSSBB method.* Let $\{b_n\}$ be a sequence of positive numbers such that

(4.2) $$b_n^{-1} + b_n/\lambda_n \to 0 \quad \text{as } n \to \infty.$$

For simplicity of exposition, in this section we also suppose that $r_n \equiv \lambda_n/b_n$ is an even integer for all $n \geq 1$. Here $b_n$ determines the volumes of the bootstrap blocks. Condition (4.2) says that the volumes of the DSSBB blocks grow to infinity, but at a rate slower than the volume of the sampling region $R_n$. The basic idea behind the formulation of the DSSBB method is to partition the sampling region $R_n$ into $d$-dimensional cubes of volume $b_n^d$ and define a version of the spatial process on each such subregion by sampling from a suitable collection of data-site-shifted blocks. Specifically, let

(4.3) $$R_n = \bigcup_{\mathbf{k} \in \mathcal{K}_n} R_n(\mathbf{k})$$

be the partition of $R_n$, where $R_n(\mathbf{k}) \equiv R_n \cap [(\mathbf{k} + [0,1)^d)b_n]$, $\mathbf{k} \in \mathbb{Z}^d$, and $\mathcal{K}_n = \{\mathbf{k} \in \mathbb{Z}^d : R_n(\mathbf{k}) \neq \varnothing\}$. For $R_n$ of (4.1), it is easy to check that $\mathcal{K}_n = \{(k_1,\ldots,k_d)' \in \mathbb{Z}^d : -r_n \leq 2k_i < r_n \text{ for } i=1,\ldots,d\}$. Next, let $\mathcal{B}_n^{[D]}(i) = \mathbf{s}_i + [0,1)^d b_n$, $1 \leq i \leq n$, denote the cubes of volume $b_n^d$ with their "lower left" end-points at the sampling sites $\mathbf{s}_1,\ldots,\mathbf{s}_n$ and define

$$\mathcal{I}_n^{[D]} = \{i : 1 \leq i \leq n, \mathbf{s}_i + [0,1)^d b_n \subset R_n\},$$

the index set of all data-site-shifted (observed) blocks of volume $b_n^d$ that are contained in the sampling region $R_n$ [cf. Figure 1(b)(i)]. Let $\{\mathcal{B}_n^{**}(\mathbf{k}) : \mathbf{k} \in \mathcal{K}_n\}$ be a random sample drawn with replacement from the collection $\{\mathcal{B}_n^{[D]}(i) : i \in \mathcal{I}_n^{[D]}\}$. Then the DSSBB version of the $Y$-process over $R_n$ is given by concatenating the observations in the resampled blocks $\{\mathcal{B}_n^{**}(\mathbf{k}) : \mathbf{k} \in \mathcal{K}_n\}$. More precisely, let $\mathcal{Y}_n(A) = \{Y(\mathbf{s}_i) : 1 \leq i \leq n, \mathbf{s}_i \in A\}$ denote the set of observations over a set of $A \subset \mathbb{R}^d$. Then, in this notation,

(4.4) $$\mathcal{Y}_n(R_n) = \{Y(\mathbf{s}_i) : i = 1,\ldots,n\},$$

the collection of all observations. The DSSBB version of $\mathcal{Y}_n(R_n)$ is defined as

(4.5) $$\mathcal{Y}_n^{**}(R_n) \equiv \bigcup_{\mathbf{k} \in \mathcal{K}_n} \mathcal{Y}_n(\mathcal{B}_n^{**}(\mathbf{k})).$$

The DSSBB version of a statistic $T_n \equiv t_n(\mathcal{Y}_n(R_n))$ is given by replacing the observations $\{Y(\mathbf{s}_i) : i = 1,\ldots,n\}$ by the resampled values $\mathcal{Y}_n^{**}(R_n)$.

When the sampling region is of a general shape, not all subregions $R_n(\mathbf{k})$ in (4.3) are necessarily cubic [cf. Figure 1(a)]. In this case, we may only use a part of a resampled cubic block that is congruent to a given subregion [cf. Figure 1(b)(ii)] and combine all the resampled observations to define the DSSBB version of the $Y$-process on all of $R_n$. (See the description of the GBBB method in Section 5 for more details.) Politis and Sherman [34] establish validity of the DSSBB method for marked point processes, where the sampling sites are generated by a collection of weakly dependent random vectors, but with the *uniform* distribution on the sampling region. In the next section we give an example in the two-dimensional case to show that if the design density is *nonuniform*, the DSSBB method fails to provide a consistent variance estimator even for the simple case where $T_n$ is the sample mean.

4.2. *Inconsistency of the DSSBB method under nonuniformity.* Suppose that $\{Y(\mathbf{s}) : \mathbf{s} \in \mathbb{R}^d\}$ is a stationary r.f. on the plane (i.e., $d=2$) and that the prototype set $R_0$ is given by $R_0 = (-1/2, 1/2)^2$ [cf. (4.1)]. Further, suppose that the sample size $n$ and the sequence $\{\lambda_n\}$ satisfy the relation

$$n = (\lambda_n^d)^2. \tag{4.6}$$

Note that (4.6) implies a mixed increasing domain spatial asymptotic structure with a nontrivial infill component (which corresponds to the case "$c = \infty$" in Theorem 1). Next, suppose that the spatial sampling density $f$ is given by

$$f_a(x_1, x_2) = g_a(x_1) \mathbb{1}_{(-1/2, 1/2)}(x_2), \tag{4.7}$$

for $a \in (4, \infty)$ (to be specified later), where $g_a(x_1)$ is a symmetric function that, on the interval $(0, \infty)$, is given by

$$g_a(x_1) = \begin{cases} a/4, & 0 < x_1 < a^{-1}, \\ \text{linear}, & a^{-1} < x_1 < 2a^{-1}, \\ a/[4(a-3)], & 2a^{-1} < x_1 < 1/2, \\ 0, & x_1 > 1/2. \end{cases} \tag{4.8}$$

For $a$ large, the sampling design puts one-half of its mass on the thin strip $(-a^{-1}, a^{-1})(-1/2, 1/2)$ and has the other one-half mass on the rest of the unit square. Next let

$$T_n = n^{-1} \sum_{i=1}^{n} Y(\mathbf{s}_i)$$

denote the sample mean and let $\mathcal{Y}_n^{**}(R_n)$ denote the DSSBB version of the observed $Y$-process based on blocks of size $b_n$. Then, following the description of Section 4.1, the DSSBB version of $T_n$ is given by

$$T_n^{**} = \text{average of the resampled observations in } \mathcal{Y}_n^{**}(R_n). \tag{4.9}$$

Suppose that we wish to estimate the variance of $T_n$ using the DSSBB method. In Proposition 1 below, we show that the scaled variance of $T_n$,

$$\check{\sigma}_n^2 \equiv \lambda_n^d \operatorname{Var}_{\cdot | \mathbf{X}}(T_n), \tag{4.10}$$

has a nonrandom limit $\sigma_\infty^2$ (say) for almost all realizations of the spatial stochastic design vectors $\mathbf{X}_1, \mathbf{X}_2, \ldots$, where $\mathrm{Var}_{\cdot|\mathbf{X}}$ denotes the conditional variance given $\mathcal{X}$. Define the DSSBB estimator of the scaled variance $\check{\sigma}_n^2$ of $T_n$ or of the limiting parameter $\sigma_\infty^2$ by

$$\hat{\sigma}_n^2 \equiv \lambda_n^d \mathrm{Var}_*(T_n^{**}),$$

where $T_n^{**}$ is the DSSBB version of $T_n$ given by (4.9) and where $\mathrm{Var}_*$ denotes the conditional variance given $\{Y(\mathbf{s}) : \mathbf{s} \in \mathbb{R}^2\} \cup \{\mathbf{X}_i : i \geq 1\}$. The following result gives the large sample behavior of $\hat{\sigma}_n^2$. Recall that in this paper, $(\Omega, \mathcal{F}, P)$ denotes the underlying probability space and hence, $P$ denotes the unconditional probability measure.

PROPOSITION 1. *Suppose that* $\{Y(\mathbf{s}) : \mathbf{s} \in \mathbb{R}^2\}$ *is a zero-mean bounded stationary r.f. with strong-mixing coefficient* $\alpha_Y(\cdot; \cdot)$ [*cf.* (3.4)] *satisfying*

(4.11) $\qquad \alpha_Y(a; b) \leq C_1 \exp(-C_2 a) b^{C_3}, \qquad a \geq 1, b \geq 1,$

*for some constants* $C_1, C_2, C_3 \in (0, \infty)$. *Also suppose that relation* (4.6) *holds and the spatial sampling density is given by* (4.7). *Then:*

(i) $\check{\sigma}_n^2 \to \sigma_\infty^2$ *a.s.* $(P_\mathbf{X})$, *where* $\sigma_\infty^2 = \int \sigma(\mathbf{s}) \, d\mathbf{s} \int f_a^2(\mathbf{s}) \, d\mathbf{s}$ *and* $\sigma(\mathbf{s}) = \mathrm{Cov}(Z(\mathbf{s}), Z(\mathbf{0}))$, $\mathbf{s} \in \mathbb{R}^2$.

(ii) *Suppose, in addition, that* $\sigma(\mathbf{0}) \neq 0$ *and* $\sigma(\mathbf{s}) \geq 0$ *for all* $\mathbf{s} \in \mathbb{R}^2$. *Then there exist* $a \in (4, \infty)$ *and* $\eta = \eta(a) \in (0, 1)$ *such that*

(4.12) $\qquad \limsup_{n \to \infty} P(\hat{\sigma}_n^2 > \eta \sigma_\infty^2) < 1.$

A proof of Proposition 1 is given in Section 9. Note that (4.12) shows that $\hat{\sigma}_n^2$ is not a consistent estimator of $\sigma_\infty^2$. Thus, the DSSBB method fails to provide a valid approximation even for very smooth functionals of the sampling distribution (e.g., the variance) of the sample mean under nonuniform spatial stochastic designs. It can be shown that a similar inconsistency result continues to hold for the pure increasing domain spatial asymptotic structure,

$$n/\lambda_n^2 \to c \in (0, \infty),$$

and also for other cases of mixed increasing domain asymptotic structure where $n/\lambda_n^2 \to \infty$ possibly at a different rate. The failure of the DSSBB method seems to be an artifact of the *interaction* between the nonuniform design density and of the additional randomness in the data-site-shifted blocks *induced by* the random locations of the sampling sites. In the next section we describe the GBBB method, which untangles the interaction of these two factors by introducing an extraneous grid to define the blocks. We also show that the GBBB method produces valid approximations for *nonuniform* sampling densities under both spatial asymptotic structures and under the *full generality* of the regression model (3.1).

## 5. The GBBB method and its consistency properties.

### 5.1. *The GBBB method.*

5.1.1. *Description of the GBBB method for a stationary random field.*
For simplicity of exposition, in this section we suppose that $\{Y(\mathbf{s}) : \mathbf{s} \in \mathbb{R}^d\}$ is a *stationary* r.f. [i.e., $\beta = 0$ in model (3.1)]. Let $\{b_n\}$ be a sequence of positive numbers satisfying (4.2), that is, $b_n \to \infty$ as $n \to \infty$ but $b_n/\lambda_n \to 0$ as $n \to \infty$. However, unlike in Section 4.1, here we do not require $\lambda_n/b_n$ to be an integer. Let $\mathcal{K}_n = \{\mathbf{k} \in \mathbb{Z}^d : \mathbf{k}b_n + [0,1)^d b_n \cap R_n \neq \varnothing\}$ denote the minimal set of indexes $\mathbf{k} \in \mathbb{Z}^d$ such that the collection $\{\mathbf{k}b_n + [0,1)^d b_n : \mathbf{k} \in \mathcal{K}_n\}$ gives a covering of $R_n$ by disjoint (hyper-)cubes of sides $b_n$. For $\mathbf{k} \in \mathcal{K}_n$, let $R_n(\mathbf{k}) = R_n \cap \{\mathbf{k}b_n + [0,1)^d b_n\}$ denote the part of $R_n$ covered by $\mathbf{k}b_n + [0,1)^d b_n$. Then $R_n$ can be expressed as a disjoint union of the sets $R_n(\mathbf{k}), \mathbf{k} \in \mathcal{K}_n$, that is,

$$(5.1) \qquad R_n = \bigcup_{\mathbf{k} \in \mathcal{K}_n} R_n(\mathbf{k})$$

[cf. (4.3)]. Note that for a nonrectangular sampling region $R_n$, the subregions $R_n(\mathbf{k})$'s lying near the boundary of $R_n$ need not be rectangular, and in general, may have different shapes. The proposed spatial block bootstrap method defines a bootstrap version of the $Y$-process on each of the subregions $R_n(\mathbf{k})$. To that end, let $\mathcal{I}_n = \{\mathbf{i} \in \mathbb{Z}^d : \mathbf{i} + [0,1)^d b_n \subset R_n\}$ denote the index set of all *overlapping* hypercubes of the form $\mathbf{i} + [0,1)^d b_n$ that are contained in the sampling region $R_n$. Also, let $\{I_\mathbf{k} : \mathbf{k} \in \mathcal{K}_n\}$ denote a collection of i.i.d. random variables with the discrete uniform distribution on the set $\mathcal{I}_n$, that is, for each $\mathbf{k} \in \mathcal{K}_n$,

$$(5.2) \qquad P_*(I_\mathbf{k} = \mathbf{i}) = \frac{1}{|\mathcal{I}_n|}, \qquad \mathbf{i} \in \mathcal{I}_n,$$

where, recall that, for any set $A$, $|A|$ denotes its size and $P_*$ denotes the conditional probability given $\sigma\langle\{Z(\mathbf{s}) : \mathbf{s} \in \mathbb{R}^d\} \cup \{\mathbf{X}_i : i \geq 1\}\rangle \equiv \mathcal{G}$. As a first step, we select the cubic blocks $I_\mathbf{k} + [0,1)^d b_n$, $\mathbf{k} \in \mathcal{K}_n$, using the bootstrap variables $\{I_\mathbf{k} : \mathbf{k} \in \mathcal{K}_n\}$. Since each subregion $R_n(\mathbf{k})$ in the partition of $R_n$ in (5.1) is contained in a cube of volume $b_n^d$, it is possible to inscribe a "copy" of $R_n(\mathbf{k})$ in the resampled block $I_\mathbf{k} + [0,1)^d b_n$ for every $\mathbf{i} \in \mathcal{I}_n$. More precisely, we define the $\mathbf{k}$th resampled block by

$$(5.3) \qquad B_n^*(\mathbf{k}) = R_n(\mathbf{k}) - \mathbf{k}b_n + I_\mathbf{k}, \qquad \mathbf{k} \in \mathcal{K}_n.$$

Note that as $R_n(\mathbf{k}) \subset \mathbf{k}b_n + [0,1)^d b_n$, $R_n(\mathbf{k}) - \mathbf{k}b_n + I_\mathbf{k} \subset I_\mathbf{k} + [0,1)^d b_n$. Further, $B_n^*(\mathbf{k})$ is congruent to $R_n(\mathbf{k})$, since they differ only by a *translation*. Next let $\mathcal{S}_n \equiv \{\mathbf{s}_1, \ldots, \mathbf{s}_n\}$ denote the collection of data-sites and for any set $A \subset \mathbb{R}^d$, let $\mathcal{Y}_n(A) \equiv \{Y(\mathbf{s}) : \mathbf{s} \in A \cap \mathcal{S}_n\}$ denote the collection of observations $Y(\mathbf{s}_i)$ over the set $A$. Then, like the DSSBB method, the GBBB version of the $Y$-process on $R_n$ is given by pasting together the observations in the resampled blocks $\{B_n^*(\mathbf{k}) : \mathbf{k} \in \mathcal{K}_n\}$ [cf. (4.5)]:

$$\mathcal{Y}_n^*(R_n) = \bigcup_{\mathbf{k} \in \mathcal{K}_n} \mathcal{Y}(B_n^*(\mathbf{k})).$$

The GBBB version of a statistic $T_n$ based on $\mathcal{Y}(R_n)$ is given by replacing variables in $\mathcal{Y}(R_n)$ with the resampled observations in $\mathcal{Y}_n^*(R_n)$. For example, for $T_n = n^{-1} \sum_{i=1}^n Y(\mathbf{s}_i)$, the sample mean, its GBBB version is given by $T_n^* =$ the mean of the resampled values in $\mathcal{Y}_n^*(R_n)$.

5.1.2. *Description of the GBBB method for the spatial regression model.* Extension of the GBBB method from the case of a stationary r.f. to the regression model (3.1) is done as follows. Let $\{B_n^*(\mathbf{k}) : \mathbf{k} \in \mathcal{K}_n\}$ denote the resampled blocks of the GBBB method, as given by (5.3). To define the bootstrap version of the $Y$-process, we first define the bootstrap version of the error process $Z(\cdot)$ on $R_n$ and then use the structural relation of the regression model (3.1) to construct the GBBB version of the $Y$-process. To that end, let

$$\hat{Z}_n(\mathbf{s}_i) = Y(\mathbf{s}_i) - \mathbf{X}(\mathbf{s}_i)' \hat{\beta}_n, \qquad i = 1, \ldots, n,$$

denote the residuals, and for $A \subset \mathbb{R}^d$, let $\hat{\mathcal{Z}}_n(A) \equiv \{\hat{Z}_n(\mathbf{s}) : \mathbf{s} \in A \cap \mathcal{S}_n\}$ denote the collection of the residuals $\hat{Z}_n(\mathbf{s}_i)$ over the set $A$, where recall that $\mathcal{S}_n = \{\mathbf{s}_1, \ldots, \mathbf{s}_n\}$. With this, we define the bootstrap version of the error process $Z(\cdot)$ over the subregion $R_n(\mathbf{k})$ by

(5.4) $$\mathcal{Z}_n^*(R_n(\mathbf{k})) \equiv \hat{\mathcal{Z}}_n(B_n^*(\mathbf{k})), \qquad \mathbf{k} \in \mathcal{K}_n.$$

As in the stationary case (cf. Section 5.1.1), the bootstrap version $\mathcal{Z}_n^*(R_n)$ of the error variables over the entire sampling region $R_n$ is now obtained by pasting together the bootstrap versions $\mathcal{Z}_n^*(R_n(\mathbf{k}))$ for all $\mathbf{k} \in \mathcal{K}_n$, that is, by

(5.5) $$\mathcal{Z}_n^*(R_n) = \bigcup_{\mathbf{k} \in \mathcal{K}_n} \mathcal{Z}_n^*(R_n(\mathbf{k})).$$

Next note that as the sampling sites $\mathbf{s}_1, \ldots, \mathbf{s}_n$ are irregularly spaced, the numbers of sampling sites in the blocks $\mathbf{i} + [0,1)^d b_n$ and $\mathbf{j} + [0,1)^d b_n$ for any two distinct indexes $\mathbf{i}, \mathbf{j} \in \mathcal{I}_n$ are typically different. Let $L_{\mathbf{k}}^* \equiv |B_n^*(\mathbf{k}) \cap \mathcal{S}_n|$ denote the number of observations in the resampled block $B_n^*(\mathbf{k})$. Then the bootstrap sample size under the present resampling scheme is given by

(5.6) $$N^* \equiv \sum_{\mathbf{k} \in \mathcal{K}_n} L_{\mathbf{k}}^*,$$

which is *random* and is, in general, different from the original sample size $n$. However, it can be shown that the expected value of $N^*$ is asymptotically equivalent to $n$. Let $\{\mathbf{s}_i^* : i = 1, \ldots, N^*\}$ denote the collection of sampling sites in the resampled blocks $\{B_n^*(\mathbf{k}) : \mathbf{k} \in \mathcal{K}_n\}$, that is,

$$(5.7) \qquad \{\mathbf{s}_i^* : i = 1, \ldots, N^*\} = \bigcup_{\mathbf{k} \in \mathcal{K}_n} [B_n^*(\mathbf{k}) \cap \mathcal{S}_n],$$

and let $\{Z_n^*(\mathbf{s}_i^*) : i = 1, \ldots, N^*\}$ denote the corresponding listing of the bootstrap error variables in the collection $\mathcal{Z}_n^*(R_n)$. Then we define the bootstrap "observations" $Y^*(\cdot)$'s by

$$(5.8) \qquad Y^*(\mathbf{s}_i^*) = \mathbf{w}(\mathbf{s}_i^*)'\hat{\beta}_n + Z_n^*(\mathbf{s}_i^*), \qquad i = 1, \ldots, N^*.$$

Given a random variable $T_n \equiv t_n(\mathcal{D}_n(R_n); \beta)$, where $\mathcal{D}_n(R_n) = \{(\mathbf{w}_n(\mathbf{s}_i)', Y(\mathbf{s}_i))' : i = 1, \ldots, n\}$ stands for the data at hand, the (overlapping) GBBB version $T_n^*$ of $T_n$ is given by replacing $\mathcal{D}_n(R_n)$ with $\{(\mathbf{w}_n(\mathbf{s}_i^*)', Y^*(\mathbf{s}_i^*))' : i = 1, \ldots, N^*\} \equiv \mathcal{D}_n^*(R_n)$ and $t_n(\cdot)$ with $t_{N^*}(\cdot)$ in the definition of $T_n$:

$$(5.9) \qquad T_n^* = t_{N^*}(\mathcal{D}_n^*(R_n); \hat{\beta}_n).$$

Next we apply the above block bootstrap method to a normalized version of the $M$-estimator:

$$(5.10) \qquad T_{1n} = \Lambda_{1n}(\hat{\beta}_n - \beta),$$

where $\Lambda_{1n} \equiv \lambda_n^{d/2} \Lambda_n'$ and $\Lambda_n$ is as in condition (C.1). Existing literature on bootstrapping $M$-estimators of regression parameters with independent errors and/or time-series errors suggests that there can be more than one way of defining the bootstrap version of the normalized $M$-estimator $T_{1n}$ (cf. [10, 18]). Here we follow an approach initially put forward by Shorack [38] in the context of bootstrapping $M$-estimators in regression models with i.i.d. errors. To define the bootstrap version of $T_{1n}$, define the variables $S_n^*(\mathbf{k}; \mathbf{t}) = \sum_{i=1}^{N^*} \mathbf{w}(\mathbf{s}_i^*) \psi(Y^*(\mathbf{s}_i^*) - \mathbf{w}(\mathbf{s}_i^*)'\mathbf{t}) \mathbb{1}(\mathbf{s}_i^* \in B_n^*(\mathbf{k}))$, $\mathbf{t} \in \mathbb{R}^p, \mathbf{k} \in \mathcal{K}_n$. Also, let

$$\hat{c}_n(\mathbf{k}) \equiv E_* \left[ \sum_{i=1}^{N^*} \mathbf{w}(\mathbf{s}_i^*) \psi(\hat{Z}_n^*(\mathbf{s}_i^*)) \mathbb{1}(\mathbf{s}_i^* \in B_n^*(\mathbf{k})) \right],$$

$\mathbf{k} \in \mathcal{K}_n$, where $E_*$ denotes the conditional expectation given $\mathcal{G}$. Then we define the bootstrap version $\beta_n^*$ of $\hat{\beta}_n$ as a measurable solution to the equation (in $\mathbf{t} \in \mathbb{R}^p$)

$$(5.11) \qquad \sum_{\mathbf{k} \in \mathcal{K}_n} [S_n^*(\mathbf{k}; \mathbf{t}) - \hat{c}_n(\mathbf{k})] = 0.$$

The motivation behind this definition is the following. For $\mathbf{t} \in \mathbb{R}^p$, $\mathbf{k} \in \mathcal{K}_n$, setting $S_n(\mathbf{k}; \mathbf{t}) \equiv \sum_{i=1}^n [\mathbf{w}(\mathbf{s}_i) \psi(Y(\mathbf{s}_i) - \mathbf{w}(\mathbf{s}_i)'\mathbf{t}) \mathbb{1}(\mathbf{s}_i \in R_n(\mathbf{k}))]$, we may rewrite the estimating equation (3.3) (defining the $M$-estimator $\hat{\beta}_n$) as

$$(5.12) \qquad \sum_{\mathbf{k} \in \mathcal{K}_n} S_n(\mathbf{k}; \mathbf{t}) = 0.$$

By (3.2), the expected value of the sum on the left-hand side of (5.12) at the true parameter value $\mathbf{t} = \beta$ is zero. However, the bootstrap versions $S_n^*(\mathbf{k}; \mathbf{t})$ of $S_n(\mathbf{k}; \mathbf{t})$ do not automatically satisfy this model restriction (in the bootstrap world). Here we introduced the centering values $\hat{c}_n(\mathbf{k})$ in the bootstrap estimating equation (5.11) precisely to ensure that the analog of this unbiasedness condition holds at the bootstrap "true" value $\hat{\beta}_n$.

With the bootstrap version $\beta_n^*$ of the $M$-estimator $\hat{\beta}_n$ defined by (5.11), we now define the bootstrap version of the normalized $M$-estimator $T_{1n}$ as

$$(5.13) \qquad T_{1n}^* = \Lambda_{1n}(\beta_n^* - \hat{\beta}_n).$$

5.1.3. *An example.* For an illustration of the main steps leading to the GBBB version of $T_{1n}^*$ of the normalized $M$-estimator $T_{1n}$, consider model (3.1), with the weight function $\mathbf{w}(\mathbf{s}) \equiv 1, \mathbf{s} \in R_n$, and $\psi(x) = x$, $x \in \mathbb{R}$. Then $\{Y(\mathbf{s}) : \mathbf{s} \in \mathbb{R}^d\}$ is a stationary r.f. with $EY(\mathbf{0}) = \mu$. In this case, $\hat{\beta}_n = \bar{Y}_n$ corresponds to the $M$-estimator of the mean $\beta \equiv \mu$, and $T_{1n}$ to the normalized sample mean

$$T_{1n} = \Lambda_{1n}(\bar{Y}_n - \mu),$$

where $\bar{Y}_n \equiv n^{-1} \sum_{i=1}^n Y(\mathbf{s}_i)$. The residuals are given by $\hat{Z}_n(\mathbf{s}_i) = Y(\mathbf{s}_i) - \bar{Y}_n$, $i = 1, \ldots, n$. Let

$$B_n(\mathbf{i}; \mathbf{k}) \equiv R_n(\mathbf{k}) - \mathbf{k}b_n + \mathbf{i}, \mathbf{i} \in \mathcal{I}_n$$

denote the blocks of "type $\mathbf{k}$" [i.e., blocks *congruent* to $R_n(\mathbf{k})$] and let $\hat{\mu}_n(\mathbf{k}) = E_*[\sum_{i=1}^n Y(\mathbf{s}_i) \mathbb{1}(\mathbf{s}_i \in B_n^*(\mathbf{k}))]$, $\mathbf{k} \in \mathcal{K}_n$. Thus, for each $\mathbf{k} \in \mathcal{K}_n$, $\hat{\mu}_n(\mathbf{k})$ is the average of all "type $\mathbf{k}$" block sums, where each sum extends over the sampling sites lying in the $\mathbf{i}$th "type $\mathbf{k}$" block $B_n(\mathbf{i}; \mathbf{k})$, $\mathbf{i} \in \mathcal{I}_n$. The centering constants $\hat{c}_n(\mathbf{k})$ are now given by $\hat{c}_n(\mathbf{k}) = \hat{\mu}_n(\mathbf{k}) - \bar{Y}_n E_* L_{\mathbf{k}}^*$, $\mathbf{k} \in \mathcal{K}_n$, where recall that $L_{\mathbf{k}}^*$ denotes the number of observations in the resampled block $B_n^*(\mathbf{k})$. With this, (5.11) can be rewritten as

$$\left[\sum_{i=1}^{N^*}(Y^*(\mathbf{s}_i^*) - t)\right] - \left[\sum_{\mathbf{k} \in \mathcal{K}_n} \hat{c}(\mathbf{k})\right] = 0.$$

Hence, the bootstrap version of the $M$-estimator is given by

$$\beta_n^* = \bar{Y}_n^* - N^{*-1} \sum_{\mathbf{k} \in \mathcal{K}_n} \hat{c}(\mathbf{k}),$$

where $\bar{Y}_n^* = \sum_{i=1}^{N^*} Y^*(\mathbf{s}_i^*)/N^*$ denotes the bootstrap sample mean. Further, by (5.13), the bootstrap version of $T_{1n}$ is given by

$$(5.14) \quad \begin{aligned} T_{1n}^* &= \Lambda_{1n}(\beta_n^* - \hat{\beta}_n) \\ &= \Lambda_{1n}\left(\bar{Y}_n^* - N^{*-1} \sum_{\mathbf{k} \in \mathcal{K}_n} \hat{c}(\mathbf{k}) - \bar{Y}_n\right) \\ &= \Lambda_{1n}(\bar{Y}_n^* - \tilde{\mu}_n), \end{aligned}$$

where $\tilde{\mu}_n = \sum_{\mathbf{k} \in \mathcal{K}_n} \hat{\mu}_n(\mathbf{k})/N^*$.

### 5.1.4. *Variants of the GBBB method.*

REMARK 5.1. Let $\mathcal{K}_{1n} = \{\mathbf{k} \in \mathcal{K}_n : \mathbf{k}b_n + [0,1)^d b_n \subset R_n\}$ denote the index set of all cubes $\mathbf{k}b_n + [0,1)^d b_n$ that are completely contained in $R_n$ and let $\mathcal{K}_{2n} = \mathcal{K}_n \setminus \mathcal{K}_{1n}$ denote the index set of all boundary hypercubes. Then $R_n(\mathbf{k}) = \mathbf{k}b_n + [0,1)^d b_n$ for all $\mathbf{k} \in \mathcal{K}_{1n}$, while for $\mathbf{k} \in \mathcal{K}_{2n}$, each $R_n(\mathbf{k})$ may have a different shape. As a consequence, the $B_n^*(\mathbf{k})$'s are of cubic shape for all $\mathbf{k} \in \mathcal{K}_{1n}$. However, for the boundary regions $R_n(\mathbf{k})$, $\mathbf{k} \in \mathcal{K}_{2n}$, the collection $\{B_n^*(\mathbf{k}) : \mathbf{k} \in \mathcal{K}_{2n}\}$ of resampled blocks may have different shapes, depending on the shapes of $R_n(\mathbf{k})$, $\mathbf{k} \in \mathcal{K}_{2n}$.

REMARK 5.2. The observation in the previous remark suggests that we may define an alternative and *simpler* version of the GBBB method, by restricting attention only to the cubic blocks while resampling. Let $\mathcal{I}_n = \{\mathbf{i} \in \mathbb{Z}^d : \mathbf{i} + [0,1)^d b_n \subset R_n\}$ be as defined earlier and let the $R_n(\mathbf{k})$'s be as in (5.1). Then, with $B_n(\mathbf{i}; \mathbf{0}) \equiv \mathbf{i} + [0,1)^d b_n$, $\{B_n(\mathbf{i}; \mathbf{0}) : \mathbf{i} \in \mathcal{I}_n\}$ is a collection of overlapping cubic blocks of sides $b_n$. For the simpler version of the GBBB, we make use of the variable $\{I_{\mathbf{k}} : \mathbf{k} \in \mathcal{K}_n\}$ of (5.2), but now select the random sample of cubic blocks $\{B_n(I_{\mathbf{k}}; \mathbf{0}) : \mathbf{k} \in \mathcal{K}_n\}$, all having the same (cubic) shape. The cubic-block GBBB (call it GBBB-CB) version of the error process $Z(\cdot)$ is now defined by concatenating the residuals on the resampled blocks:

$$\tilde{\mathcal{Z}}^*(\tilde{R}_n) = \bigcup_{\mathbf{k} \in \mathcal{K}_n} \hat{\mathcal{Z}}(B_n(I_{\mathbf{k}}; \mathbf{0})),$$

where $\tilde{R}_n = \bigcup_{\mathbf{k} \in \mathcal{K}_n} \{\mathbf{k} + [0,1)^d b_n\}$ is a covering of the sampling region $R_n$ by cubes of sides $b_n$ as its building blocks. The GBBB-CB version of a random variable $T_n \equiv t_n(\mathcal{D}_n(R_n); \beta)$ is given by

$$\tilde{T}_n^* = t_{\tilde{N}^*}(\tilde{\mathcal{D}}_n^*(\tilde{R}_n); \hat{\beta}_n),$$

where $\{\tilde{\mathbf{s}}_i^* : i = 1, \ldots, \tilde{N}^*\}$ is the set of data-sites in the resampled blocks $B_n(I_{\mathbf{k}}; \mathbf{0})$, $\mathbf{k} \in \mathcal{K}_n$, and $\tilde{\mathcal{D}}_n^*(\tilde{R}_n) = \{(\mathbf{w}(\tilde{\mathbf{s}}_i^*)', \tilde{Y}^*(\tilde{\mathbf{s}}_i^*))' : i = 1, \ldots, \tilde{N}^*\}$ is the GBBB-CB version of $\mathcal{D}(R_n)$, and where $\{\tilde{Y}(\tilde{\mathbf{s}}_i^*) : i = 1, \ldots, \tilde{N}^*\} \equiv \tilde{\mathcal{Y}}_n^*(\tilde{R}_n)$ is the collection of GBBB-CB observations defined using the resampled error variables in $\tilde{\mathcal{Z}}^*(\tilde{R}_n)$ as in (5.8). Thus the GBBB-CB version of $\hat{\beta}_n$ is defined as a solution $\tilde{\beta}_n^*$ of the equation [cf. (5.11)]

$$\sum_{\mathbf{k} \in \mathcal{K}_n} [\tilde{S}_n^*(\mathbf{k}; \mathbf{t}) - \tilde{c}_n(\mathbf{k})] = 0,$$

where $\mathbf{t} \in \mathbb{R}^p$, $\tilde{S}_n^*(\mathbf{k}; \mathbf{t})$ is the sum of all $[\mathbf{w}(\tilde{\mathbf{s}}_i^*)\psi(\tilde{\mathcal{Y}}^*(\tilde{\mathbf{s}}_i^*) - \mathbf{w}(\tilde{\mathbf{s}}_i^*)'\mathbf{t})]$ over the data-sites $\tilde{\mathbf{s}}_i^*$ lying in the $\mathbf{k}$th resampled GBBB-CB block $B_n(I_{\mathbf{k}}; \mathbf{0})$, and $\tilde{c}_n(\mathbf{k}) = |\mathcal{I}_n|^{-1} \sum_{\mathbf{i} \in \mathcal{I}_n} \sum_{j=1}^n \mathbf{w}(\mathbf{s}_j) \psi(\hat{Z}_n(\mathbf{s}_j)) \mathbb{1}(\mathbf{s}_j \in B_n(\mathbf{i}; \mathbf{0}))$, $\mathbf{k} \in \mathcal{K}_n$. Further, the GBBB-CB version of $T_{1n}$ is given by

$$\tilde{T}_{1n}^* = \Lambda_{1n}(\tilde{\beta}_n^* - \hat{\beta}_n),$$

where recall that $\Lambda_{1n} = \lambda_n^{d/2} \Lambda_n'$. It can be shown, by recasting the arguments in the proof of Theorem 2, that the GBBB-CB method provides a valid approximation to the distribution of $T_{1n}$ under the *same* set of regularity conditions as Theorem 2 below.

REMARK 5.3. It is evident from the formulation of the GBBB method that, like the time-series case (cf. [7]), we may similarly formulate a nonoverlapping version of the GBBB method as well. In this case, for each $\mathbf{k} \in \mathcal{K}_n$, we need to select a block at random, with replacement from the collection $\{\bar{B}_n(\mathbf{j}, \mathbf{k}) : \mathbf{j} \in \mathcal{J}_n\}$ of disjoint blocks of "type $\mathbf{k}$," where $\mathcal{J}_n = \{\mathbf{j} \in \mathbb{Z}^d : (\mathbf{j} + [0,1)^d) b_n \subset R_n\}$ and $\bar{B}_n(\mathbf{j}, \mathbf{k}) = R_n(\mathbf{k}) - \mathbf{k} + \mathbf{j} b_n$. Note that a similar modification can also be made to the overlapping GBBB-CB method to formulate a nonoverlapping version of the GBBB-CB method.

5.2. *Theoretical validity of the GBBB.* Without loss of generality, we shall suppose that $\{\mathbf{X}_n\}_{n \geq 1}$, $\{Z(\mathbf{s}) : \mathbf{s} \in \mathbb{R}^d\}$ and the bootstrap variables $\{I_\mathbf{k} : \mathbf{k} \in \mathcal{K}_n\}$, $n \geq 1$, are all defined on a common probability space $(\Omega, \mathcal{F}, P)$. Let $P_*$, $E_*$ and $\text{Var}_*$ denote the conditional probability, conditional expectation and conditional variance given $\mathcal{G} \equiv \sigma \langle \{\mathbf{X}_n : n \geq 1\} \cup \{Z(\mathbf{s}) : \mathbf{s} \in \mathbb{R}^d\} \rangle$. For simplicity of exposition, we suppose that the solution to (3.3) is unique. When this assumption does not hold, the result is valid for a suitable sequence of solutions of (3.11), as in Remark 3.1. The main result of this section is the following.

THEOREM 2. *Suppose that conditions* (C.1)–(C.5) *hold with $r = 3$. Also suppose that for some $\delta \in (0, 1)$,*

(5.15) $$b_n^{-1} + \lambda_n^{\delta - 1} b_n = o(1) \qquad \text{as } n \to \infty$$

*and that $m_{0n} = O(1)$ as $n \to \infty$. Then*

$$\sup_{A \in \mathcal{C}} |P_*(T_{1n}^* \in A) - P_{\cdot|\mathbf{X}}(T_{1n} \in A)| \to 0 \qquad \text{in } P_{\cdot|\mathbf{X}}\text{-probability a.s. } (P_\mathbf{X})$$

(5.16)
*under both the cases "$c \in (0, \infty)$" and "$c = \infty$," where $\mathcal{C}$ is the collection of all measurable convex sets in $\mathbb{R}^p$.*

Theorem 2 shows that under some general regularity conditions on the underlying spatial process, the GBBB method produces valid approximations to the distributions of $M$-estimators under both pure- and mixed-increasing domain asymptotic structures. Furthermore, unlike the DSSBB method, it remains valid for nonuniform spatial sampling designs. The bootstrap approximation generated by the GBBB method adapts itself to the various forms of the large sample behavior of the normalized $M$-estimators arising from different combinations of various spatial asymptotic structures, spatial sampling designs and weight functions, correlation structures of the underlying error process. As a result, the GBBB method serves as a unified, valid block bootstrap method for a wide range of irregularly spaced spatial data.

A second important aspect of Theorem 2 is its validity for a class of strongly dependent spatial data. Note that in the presence of the infill component, the variance of the sample sum grows at a faster rate than the sample size. As a result, the infill component of the spatial asymptotic structure induces conditions of long-range dependence in the data. Although block bootstrap methods are known to fail for certain classes of equi-spaced long-range dependent time-series data (cf. [19]), the results of this paper show that GBBB produces valid approximations even in the presence of long-range dependence here.

Although Theorem 2 remains valid for a wide range of block sizes [cf. (5.15)], the accuracy of the GBBB estimator depends on the particular block size employed in practice. At this point, we do not know the theoretical optimal spatial block size. Note that in addition to the geometry of the sampling region and the covariance structure of the underlying random process, the optimal block length for irregularly spaced spatial data would also depend on the sampling density and the spatial asymptotic structure. Quantifying the exact form of interaction among these various factors remains a challenging problem for future investigation. In the next section we investigate finite sample properties of the GBBB method through a simulation study and provide some guidelines for choosing reasonable block sizes in practice.

**6. Simulation results.** The factors we vary in the simulation consist of the sampling region size, the spatial sampling design, the sample size, the range of dependence, and the bootstrap block size in the GBBB method. We consider $R_n = \lambda_n R_0$, where $R_0 = (-1/2, 1/2]^2$ and the scaling factor $\lambda_n = 12$ or $24$. Two types of probability density function $f$ for the spatial sampling design are used. One is a uniform distribution over $R_0$ and the other a mixture of two bivariate normal distributions, $0.5N((0,0)', I_2) + 0.5N((1/4, 1/4)', 2I_2)$, truncated outside $R_0$, where $I_2$ is a $2 \times 2$ identity matrix. Sample sizes are chosen to be $n = 100, 400, 900$ and the sampling sites $\mathbf{s}_1, \ldots, \mathbf{s}_n$ are generated according to $f$. Now, we consider the class of spatial regression models in (3.1). For simplicity, we consider a simple linear regression $\mathbf{w}(\mathbf{s}) = (1, x(\mathbf{s}))'$ and $\beta = (\beta_0, \beta_1)'$, where the covariate $x(\cdot) = 0$ or 1. The true regression parameters are set at $\beta_0 = 25$ and $\beta_1 = -5$. The error process $Z(\cdot)$, on the other hand, follows a zero-mean Gaussian process with a spherical variogram that has a unit sill, range $r$, and no nugget effect. Two values of range of dependence are considered: $r = 2, 4$. Thus the simulated observations are generated according to $Y(\mathbf{s}) = \beta_0 + \beta_1 x(\mathbf{s}) + Z(\mathbf{s})$. Again, for simplicity, we focus on the identity score function $\psi(x) = x$. To carry out the GBBB method, we consider bootstrap block sizes $b_n = 2, 4, 6$ for $\lambda_n = 12$ and $b_n = 4, 6, 8, 12$ for $\lambda_n = 24$. We also implement an empirical rule for choosing the block size [14]. First we apply the GBBB method to obtain the variance estimate $\hat{\psi}$ of the $\beta$ estimates. Then for $I$ preselected subregions, we apply the GBBB method to obtain the variance estimate $\hat{\psi}_i$. The block size that minimizes the mean square error $I^{-1} \sum_{i=1}^{I} (\hat{\psi}_i - \hat{\psi})^2$ is denoted by $b_n^*$ and the block size selected is $b_n^*$ scaled by the root of the scaling factor $\lambda_n$ relative to the subregion size.

For each combination of the factors $\lambda_n, f, n, r$, we generate $S = 500$ samples. For each sample, we obtain the ordinary least squares estimates $\hat{\beta}_n$ and perform $M = 1000$ resamples for various $b_n$ to estimate the variance of $\hat{\beta}_n$. Based on the $S = 500$ samples, we estimate the true variance of $\hat{\beta}_n$. Thus, we can compute the root mean square errors (MSE) of the variance estimates (see Tables 1 and 2). For a given sampling region with a scaling factor $\lambda_n$, the root MSE tends to be smaller when the sample size $n$ is larger, the range of dependence $r$ is smaller, and the bootstrap block size $b_n$ is about the same as the range of dependence. As the scaling factor $\lambda_n$ increases from 12 to 24, the root MSE becomes even smaller. The simulation results are fairly similar for the two types of spatial sampling design. The empirical rule for choosing the block size selects the block size that has the smallest root MSE in most cases.

Further, for each sample we construct a 90% confidence interval for $\beta_0$ and $\beta_1$. Since we know the true values of $\beta_0$ and $\beta_1$, we can compute the nominal coverage probability, based on the coverage of $S = 500$ confidence intervals for each regression parameter (see Tables 3 and 4). For a given scaling factor $\lambda_n$, the coverage probability tends to be larger when the sample size $n$ is larger, the range of dependence $r$ is smaller, and the bootstrap block size $b_n$ is about the same as the range of dependence. As the scaling factor $\lambda_n$ increases from 12 to 24, the coverage probability becomes higher, approaching the asymptotic confidence level of 90%. The simulation results are again similar for the two types of spatial sampling design. The empirical rule for choosing the block size selects the block size that has the largest coverage probability in most cases.

For the sake of comparison, we implement the DSSBB method for the uniform spatial sampling design. Tables 1 and 3 show that the root MSEs using the DSSBB method are larger than those using the GBBB method when the scaling factor $\lambda_n$ is small and the two methods give comparable root MSEs when the scaling factor $\lambda_n$ is large, for the different sample

TABLE 1
*Estimated root MSE of the GBBB variance estimators, with sample sizes $n = 100, 400, 900$, scaling factors $\lambda_n = 12$ with block sizes $b_n = 2, 4, 6$, and $\lambda_n = 24$ with block sizes $b_n = 4, 6, 8, 12$, and a spherical variogram of ranges $2, 4, 10$, based on $S = 500$ simulations, each with $M = 1000$ bootstrap resamples*

| | | | $\lambda_n = 12$ | | | | $\lambda_n = 24$ | | | | |
| | | | GBBB | | DSSBB | | | GBBB | | DSSBB | |
| $n$ | $r$ | $b_n$ | $\beta_0$ | $\beta_1$ | $\beta_0$ | $\beta_1$ | $b_n$ | $\beta_0$ | $\beta_1$ | $\beta_0$ | $\beta_1$ |
|---|---|---|---|---|---|---|---|---|---|---|---|
| 100 | 2 | 2* | 0.011 | 0.011 | 0.018 | 0.015 | 4* | 0.008 | 0.010 | 0.009 | 0.013 |
| | | 4 | 0.015 | 0.015 | 0.019 | 0.017 | 6 | 0.009 | 0.013 | 0.009 | 0.014 |
| | | 6 | 0.019 | 0.020 | 0.022 | 0.022 | 8 | 0.011 | 0.016 | 0.011 | 0.017 |
| | | — | — | — | — | — | 12 | 0.013 | 0.018 | 0.014 | 0.021 |
| | 4 | 2* | 0.046 | 0.013 | 0.054 | 0.014 | 4* | 0.011 | 0.011 | 0.018 | 0.013 |
| | | 4 | 0.042 | 0.017 | 0.047 | 0.016 | 6 | 0.013 | 0.014 | 0.017 | 0.014 |
| | | 6 | 0.049 | 0.022 | 0.055 | 0.021 | 8 | 0.015 | 0.017 | 0.018 | 0.016 |
| | | — | — | — | — | — | 12 | 0.019 | 0.022 | 0.022 | 0.020 |
| 400 | 2 | 2* | 0.008 | 0.003 | 0.009 | 0.003 | 4* | 0.003 | 0.002 | 0.003 | 0.003 |
| | | 4 | 0.009 | 0.004 | 0.009 | 0.004 | 6 | 0.003 | 0.003 | 0.003 | 0.003 |
| | | 6 | 0.012 | 0.006 | 0.012 | 0.006 | 8 | 0.004 | 0.004 | 0.004 | 0.004 |
| | | — | — | — | — | — | 12 | 0.006 | 0.005 | 0.005 | 0.006 |
| | 4 | 2* | 0.040 | 0.003 | 0.043 | 0.003 | 4* | 0.009 | 0.003 | 0.010 | 0.003 |
| | | 4 | 0.032 | 0.004 | 0.036 | 0.004 | 6 | 0.009 | 0.003 | 0.010 | 0.003 |
| | | 6 | 0.038 | 0.005 | 0.042 | 0.005 | 8 | 0.010 | 0.004 | 0.011 | 0.004 |
| | | — | — | — | — | — | 12 | 0.013 | 0.005 | 0.013 | 0.005 |
| 900 | 2 | 2* | 0.008 | 0.001 | 0.008 | 0.001 | 4* | 0.002 | 0.001 | 0.002 | 0.001 |
| | | 4 | 0.008 | 0.002 | 0.009 | 0.002 | 6 | 0.002 | 0.001 | 0.002 | 0.001 |
| | | 6 | 0.011 | 0.002 | 0.011 | 0.002 | 8 | 0.002 | 0.002 | 0.003 | 0.002 |
| | | — | — | — | — | — | 12 | 0.003 | 0.003 | 0.004 | 0.002 |
| | 4 | 2* | 0.042 | 0.001 | 0.044 | 0.001 | 4* | 0.008 | 0.001 | 0.010 | 0.001 |
| | | 4 | 0.033 | 0.002 | 0.036 | 0.002 | 6 | 0.007 | 0.002 | 0.009 | 0.002 |
| | | 6 | 0.039 | 0.002 | 0.042 | 0.002 | 8 | 0.008 | 0.002 | 0.010 | 0.002 |
| | | — | — | — | — | — | 12 | 0.011 | 0.002 | 0.013 | 0.003 |

The design density is a *uniform distribution*. Also reported are the estimated root MSE of the DSSBB variance estimators. The block size selected most often by the Hall, Horowitz and Jing [14] empirical rule is marked by *.

sizes, ranges of dependence and bootstrap block sizes. However, the coverage probability using the DSSBB method is quite a bit lower than that using the GBBB method, except for a few cases when $\lambda_n = 24$ and $n = 400, 900$.

Finally, we show the boxplots of the $S = 500$ GBBB variance estimates in Figures 2 and 3. In general, there is less variation and skewness when the sample size $n$ is larger, the range of dependence $r$ is smaller, and the scaling factor $\lambda_n$ is larger. There is no distinct difference in the boxplots for the two types of spatial sampling design.

TABLE 2
*Estimated root MSE of the GBBB variance estimators, with sample sizes $n = 100, 400, 900$, scaling factors $\lambda_n = 12$ with block sizes $b_n = 2, 4, 6$, and $\lambda_n = 24$ with block sizes $b_n = 4, 6, 8, 12$, and a spherical variogram of ranges $2, 4, 10$, based on $S = 500$ simulations, each with $M = 1000$ bootstrap resamples*

| $n$ | $r$ | $\lambda_n = 12$ | | | $\lambda_n = 24$ | | |
|---|---|---|---|---|---|---|---|
| | | $b_n$ | $\beta_0$ | $\beta_1$ | $b_n$ | $\beta_0$ | $\beta_1$ |
| 100 | 2 | 2* | 0.014 | 0.012 | 4 | 0.010 | 0.011 |
| | | 4 | 0.018 | 0.016 | 6* | 0.012 | 0.014 |
| | | 6 | 0.022 | 0.022 | 8 | 0.014 | 0.017 |
| | | — | — | — | 12 | 0.018 | 0.022 |
| | 4 | 2 | 0.057 | 0.013 | 4 | 0.021 | 0.011 |
| | | 4* | 0.053 | 0.015 | 6* | 0.022 | 0.014 |
| | | 6 | 0.062 | 0.018 | 8 | 0.025 | 0.016 |
| | | — | — | — | 12 | 0.031 | 0.021 |
| 400 | 2 | 2* | 0.013 | 0.003 | 4 | 0.004 | 0.003 |
| | | 4 | 0.014 | 0.004 | 6* | 0.005 | 0.004 |
| | | 6 | 0.017 | 0.006 | 8 | 0.006 | 0.005 |
| | | — | — | — | 12 | 0.007 | 0.006 |
| | 4 | 2* | 0.058 | 0.003 | 4 | 0.017 | 0.003 |
| | | 4 | 0.050 | 0.004 | 6* | 0.018 | 0.004 |
| | | 6 | 0.057 | 0.005 | 8 | 0.020 | 0.004 |
| | | — | — | — | 12 | 0.024 | 0.006 |
| 900 | 2 | 2* | 0.012 | 0.001 | 4* | 0.004 | 0.001 |
| | | 4 | 0.012 | 0.002 | 6 | 0.004 | 0.002 |
| | | 6 | 0.015 | 0.002 | 8 | 0.005 | 0.002 |
| | | — | — | — | 12 | 0.006 | 0.003 |
| | 4 | 2* | 0.057 | 0.001 | 4* | 0.015 | 0.001 |
| | | 4 | 0.048 | 0.002 | 6 | 0.015 | 0.002 |
| | | 6 | 0.055 | 0.002 | 8 | 0.017 | 0.002 |
| | | — | — | — | 12 | 0.021 | 0.002 |

The design density is a *mixture of two normal distributions*. The block size selected most often by the Hall, Horowitz and Jing [14] empirical rule is marked by *.

**7. Example.** American elk are of high conservation and management interest (cf. [1]). This example concerns the relationship between elk and landscape characteristics in northern Wisconsin, where elk were reintroduced in the early 1990s. There are a total of 514 sampling locations in a 15 km × 12 km stand (Figure 4). At each location, a 30 m × 30 m plot was surveyed for various plant abundance, including the biomass of grass, forbs and shrubs. Among the 514 locations, elk were known to be absent at 250 plots and present at the remaining 264 plots. One question of interest was to compare the biomass indices at locations where elk were absent versus those

TABLE 3
*Estimated nominal coverage of the 90% GBBB confidence intervals for $\beta_0, \beta_1$ with sample sizes $n = 100, 400, 900$, scaling factors $\lambda_n = 12$ with block sizes $b_n = 2, 4, 6$, and $\lambda_n = 24$ with block sizes $b_n = 4, 6, 8, 12$, and a spherical variogram of ranges $2, 4, 10$, based on $S = 500$ simulations, each with $M = 1000$ bootstrap resamples*

|  |  |  | $\lambda_n = 12$ |  |  |  |  | $\lambda_n = 24$ |  |  |  |
|  |  |  | GBBB |  | DSSBB |  |  | GBBB |  | DSSBB |  |
| $n$ | $r$ | $b_n$ | $\beta_0$ | $\beta_1$ | $\beta_0$ | $\beta_1$ | $b_n$ | $\beta_0$ | $\beta_1$ | $\beta_0$ | $\beta_1$ |
|---|---|---|---|---|---|---|---|---|---|---|---|
| 100 | 2 | 2* | 0.876 | 0.884 | 0.782 | 0.82 | 4* | 0.88 | 0.892 | 0.782 | 0.84 |
|  |  | 4 | 0.858 | 0.85 | 0.788 | 0.794 | 6 | 0.864 | 0.888 | 0.816 | 0.854 |
|  |  | 6 | 0.774 | 0.78 | 0.704 | 0.702 | 8 | 0.84 | 0.86 | 0.794 | 0.844 |
|  |  | — | — | — | — | — | 12 | 0.736 | 0.774 | 0.686 | 0.752 |
|  | 4 | 2* | 0.704 | 0.862 | 0.58 | 0.838 | 4* | 0.85 | 0.884 | 0.72 | 0.808 |
|  |  | 4 | 0.75 | 0.832 | 0.656 | 0.826 | 6 | 0.838 | 0.862 | 0.766 | 0.826 |
|  |  | 6 | 0.664 | 0.746 | 0.556 | 0.722 | 8 | 0.818 | 0.846 | 0.758 | 0.81 |
|  |  | — | — | — | — | — | 12 | 0.722 | 0.756 | 0.668 | 0.71 |
| 400 | 2 | 2* | 0.808 | 0.896 | 0.75 | 0.838 | 4* | 0.85 | 0.902 | 0.856 | 0.878 |
|  |  | 4 | 0.814 | 0.854 | 0.762 | 0.812 | 6 | 0.83 | 0.892 | 0.85 | 0.876 |
|  |  | 6 | 0.702 | 0.746 | 0.686 | 0.714 | 8 | 0.818 | 0.862 | 0.832 | 0.832 |
|  |  | — | — | — | — | — | 12 | 0.724 | 0.754 | 0.742 | 0.728 |
|  | 4 | 2* | 0.648 | 0.896 | 0.596 | 0.86 | 4* | 0.822 | 0.872 | 0.774 | 0.866 |
|  |  | 4 | 0.756 | 0.85 | 0.726 | 0.818 | 6 | 0.836 | 0.846 | 0.804 | 0.86 |
|  |  | 6 | 0.67 | 0.764 | 0.636 | 0.716 | 8 | 0.816 | 0.814 | 0.804 | 0.826 |
|  |  | — | — | — | — | — | 12 | 0.72 | 0.734 | 0.68 | 0.73 |
| 900 | 2 | 2* | 0.788 | 0.898 | 0.802 | 0.854 | 4* | 0.856 | 0.878 | 0.848 | 0.892 |
|  |  | 4 | 0.79 | 0.874 | 0.81 | 0.828 | 6 | 0.85 | 0.858 | 0.852 | 0.872 |
|  |  | 6 | 0.688 | 0.768 | 0.694 | 0.75 | 8 | 0.818 | 0.836 | 0.816 | 0.842 |
|  |  | — | — | — | — | — | 12 | 0.726 | 0.72 | 0.724 | 0.75 |
|  | 4 | 2* | 0.594 | 0.88 | 0.566 | 0.87 | 4* | 0.79 | 0.862 | 0.758 | 0.864 |
|  |  | 4 | 0.702 | 0.82 | 0.662 | 0.81 | 6 | 0.806 | 0.84 | 0.774 | 0.84 |
|  |  | 6 | 0.662 | 0.758 | 0.584 | 0.726 | 8 | 0.79 | 0.824 | 0.772 | 0.806 |
|  |  | — | — | — | — | — | 12 | 0.702 | 0.74 | 0.674 | 0.716 |

The design density is a *uniform distribution*. Also reported are the estimated nominal coverage of the 90% DSSBB confidence intervals for $\beta_0, \beta_1$. The block size selected most often by the Hall, Horowitz and Jing [14] empirical rule is marked by *.

where elk were present. The answers to the question could help understand the forage preferences and the movement patterns of elk.

We model the biomass indices by a linear regression, $z(\mathbf{s}) = \beta_0 + \beta_1 x(\mathbf{s}) + \varepsilon(\mathbf{s})$, where $z(\mathbf{s})$ is a biomass index, $x(\mathbf{s}) = \mathbb{1}(\text{elk is present at } \mathbf{s})$ and $\varepsilon(\cdot)$ is an error process. Here we account for the spatial dependence in the error process and draw inference on $\beta$ by the spatial block bootstrap method developed above. The empirical rule shown in Section 6 is used to determine the block size. For grass biomass, the block size selected is 6 km. The 90% confidence intervals for $\beta_0$ and $\beta_1$ are $(23.33, 27.04)$ and $(-4.25, -0.35)$ and

TABLE 4
*Estimated nominal coverage of the 90% GBBB confidence intervals for $\beta_0, \beta_1$ with sample sizes $n = 100, 400, 900$, scaling factors $\lambda_n = 12$ with block sizes $b_n = 2, 4, 6$, and $\lambda_n = 24$ with block sizes $b_n = 4, 6, 8, 12$, and a spherical variogram of ranges $2, 4, 10$, based on $S = 500$ simulations, each with $M = 1000$ bootstrap resamples*

| $n$ | $r$ | $\lambda_n = 12$ | | | $\lambda_n = 24$ | | |
|---|---|---|---|---|---|---|---|
| | | $b_n$ | $\beta_0$ | $\beta_1$ | $b_n$ | $\beta_0$ | $\beta_1$ |
| 100 | 2 | 2* | 0.84 | 0.888 | 4 | 0.848 | 0.872 |
| | | 4 | 0.816 | 0.832 | 6* | 0.82 | 0.822 |
| | | 6 | 0.74 | 0.744 | 8 | 0.798 | 0.79 |
| | | — | — | — | 12 | 0.7 | 0.688 |
| | 4 | 2 | 0.648 | 0.876 | 4 | 0.746 | 0.858 |
| | | 4* | 0.676 | 0.816 | 6* | 0.73 | 0.834 |
| | | 6 | 0.562 | 0.728 | 8 | 0.698 | 0.81 |
| | | — | — | — | 12 | 0.572 | 0.698 |
| 400 | 2 | 2* | 0.776 | 0.876 | 4 | 0.856 | 0.836 |
| | | 4 | 0.76 | 0.816 | 6* | 0.818 | 0.836 |
| | | 6 | 0.66 | 0.712 | 8 | 0.792 | 0.812 |
| | | — | — | — | 12 | 0.65 | 0.712 |
| | 4 | 2* | 0.582 | 0.898 | 4 | 0.742 | 0.874 |
| | | 4 | 0.66 | 0.834 | 6* | 0.74 | 0.86 |
| | | 6 | 0.538 | 0.732 | 8 | 0.702 | 0.822 |
| | | — | — | — | 12 | 0.584 | 0.716 |
| 900 | 2 | 2* | 0.75 | 0.878 | 4* | 0.83 | 0.886 |
| | | 4 | 0.74 | 0.828 | 6 | 0.798 | 0.844 |
| | | 6 | 0.616 | 0.744 | 8 | 0.746 | 0.796 |
| | | — | — | — | 12 | 0.628 | 0.692 |
| | 4 | 2* | 0.572 | 0.872 | 4* | 0.744 | 0.882 |
| | | 4 | 0.674 | 0.804 | 6 | 0.746 | 0.848 |
| | | 6 | 0.592 | 0.702 | 8 | 0.696 | 0.822 |
| | | — | — | — | 12 | 0.582 | 0.664 |

The design density is a *mixture of two normal distributions*. The block size selected most often by the Hall, Horowitz and Jing [14] empirical rule is marked by *.

the $p$-value for testing $H_0 : \beta_1 = 0$ is 0.056. There is some evidence that the grass biomass was different where elk were present. For the forb biomass, the block size selected is 3 km. The 90% confidence intervals for $\beta_0$ and $\beta_1$ are $(14.76, 16.37)$ and $(2.09, 4.84)$ and the $p$-value for testing $H_0 : \beta_1 = 0$ is less than 0.001. There is strong evidence that the forb biomass was different where elk were present. For the shrub biomass, the block size selected was 5 km. The 90% confidence intervals for $\beta_0$ and $\beta_1$ are $(12.74, 14.77)$ and $(0.62, 2.98)$ and the $p$-value for testing $H_0 : \beta_1 = 0$ is 0.020. There is weak evidence that the shrub biomass was different where elk were present.

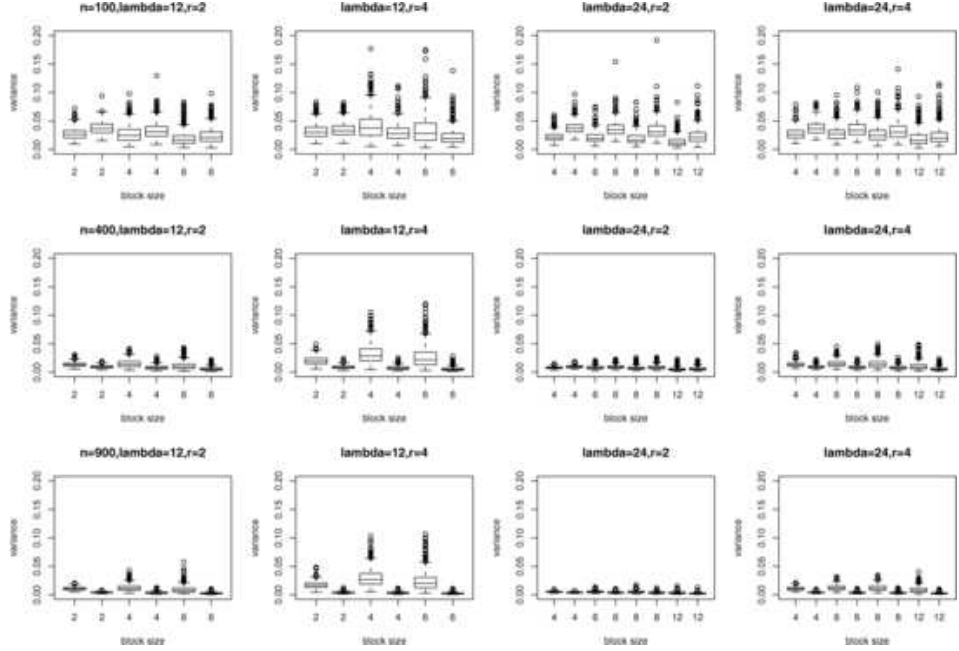

FIG. 2. *Boxplots of GBBB variance estimates for combinations of $\lambda_n = 12, 24$, $n = 100, 400, 900$, $r = 2, 4$, and various $b_n$. Two boxplots are placed next to each other for $\beta_0$ and $\beta_1$ within each subfigure. The sampling design is a uniform distribution. Simulation size is $S = 500$ and resample size is $M = 1000$.*

If standard linear regression is used [i.e., under the assumption of normal i.i.d. $\varepsilon(\cdot)$], the $p$-values for testing $H_0 : \beta_1 = 0$ are $0.030, < 0.0001$ and $0.031$, the 90% confidence intervals for $\beta_0$ are $(23.93, 26.89)$, $(14.89, 16.12)$ and $(12.71, 14.78)$, and the 90% confidence intervals for $\beta_1$ are $(-4.81, -0.67)$, $(2.21, 3.93)$ and $(0.46, 3.36)$, for grass, forb and shrub biomass, respectively. After accounting for the spatial dependence, the 90% confidence intervals could be either wider or narrower and the $p$-values could be either above or below those from the standard linear regression.

In the next three sections we give the proofs of Theorem 1 (on the asymptotic distribution of the $M$-estimators), Proposition 1 (on inconsistency of the DSSBB method) and Theorem 2 (on validity of the GBBB method).

**8. Proof of Theorem 1.** Let $\mathbf{v}(\mathbf{s}) = \Lambda_n^{-1}\mathbf{w}(\mathbf{s})$, $\mathbf{s} \in R_n$, $\eta_n^d = n^{-1}\lambda_n^d$, $C_n^{-1} = \eta_n^d \lambda_n^{-d/2} \Lambda_n^{-1}$, $n \geq 1$. Also recall that $m_{0n}^2 = \sup\{\|\mathbf{v}(\mathbf{s})\| : \mathbf{s} \in R_n\}$, $n \geq 1$. Let $\text{tr}(B)$ denote the trace of a square matrix $B$. Let $E_{\mathbf{X}}$ denote expectation with respect to the joint distribution $P_{\mathbf{X}}$ of $\mathbf{X}_1, \mathbf{X}_2, \ldots$. Let $\mathcal{U} = [0,1)^d$. Write $\mathbb{N} = \{1, 2, \ldots\}$, $\alpha_1(t) = t^{-\eta_1}$ and $g(t) = t^{\eta_2}$, $t > 0$. Set $A(r, \delta) \equiv \int_1^\infty y^{(2r-1)d-1} \alpha_1(y)^{\delta/2r+\delta} \, dy$ for $r \in \mathbb{N}, \delta \in (0, \infty)$, and let $m_n =$ $\max\{\log n, n\lambda_n^{-d}\}$, $n \geq 1$. For a nonempty set $A$ and a function $g: A \to \mathbb{R}$, let $\|g\|_\infty = \sup\{|g(x)|: x \in A\}$. Recall that $\mathcal{X} = \sigma\langle\{\mathbf{X}_1, \mathbf{X}_2, \ldots\}\rangle$ and that $\mathbb{1}(\cdot)$ denotes the indicator function and that $C, C(\cdot)$ denote generic positive constants not depending on $n$ and $\omega \in \Omega$. Furthermore, unless otherwise specified, limits in order symbols are taken letting $n \to \infty$.

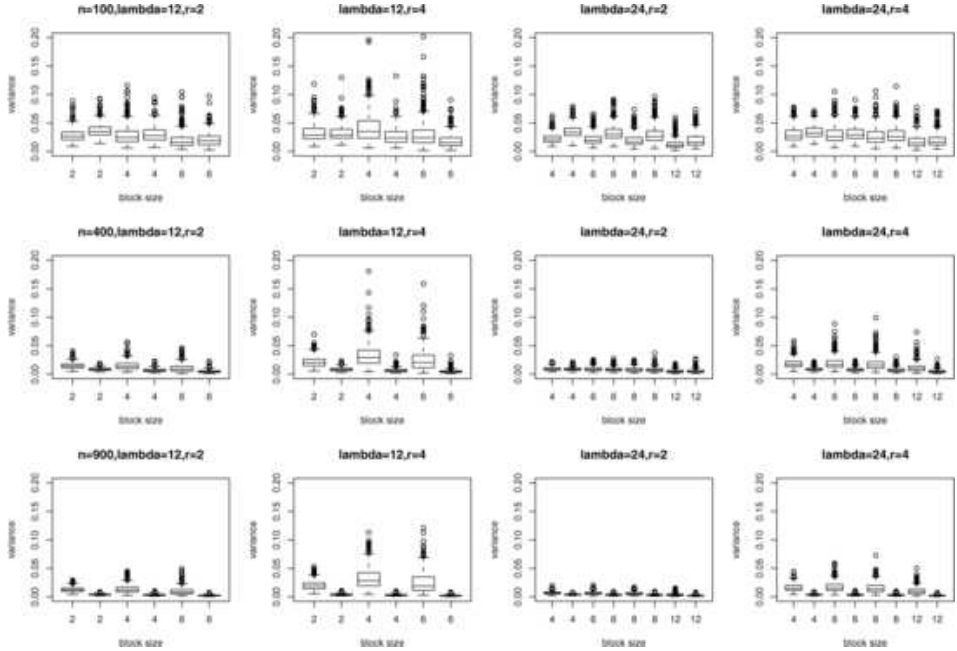

Fig. 3. *Boxplots of GBBB variance estimates for combinations of $\lambda_n = 12, 24$, $n = 100, 400, 900$, $r = 2, 4$, and various $b_n$. Two boxplots are placed next to each other for $\beta_0$ and $\beta_1$ within each subfigure. The sampling design is a mixture of two normal distributions. Simulation size is $S = 500$ and resample size is $M = 1000$.*

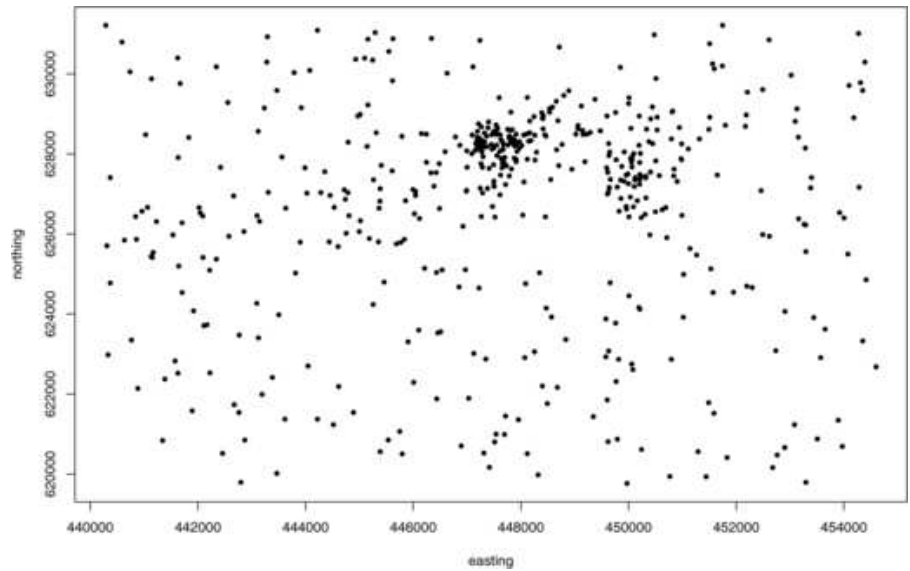

Fig. 4. *Sampling locations in a study area in northern Wisconsin. The x-axis and the y-axis are in the east–west and north–south directions (in m).*

LEMMA 1. *Let $g: \mathbb{R} \to \mathbb{R}$ be a Borel-measurable function satisfying $E|g(Z(\mathbf{0}))| < \infty$ and $Eg(Z(\mathbf{0})) = 0$. Also, let $a_{in}(\mathbf{X})$, $i = 1, \ldots, n$, be $\mathcal{X}$-measurable random variables on $(\Omega, \mathcal{F}, P)$ satisfying*

$$\text{(8.1)} \qquad \sum_{i=1}^n |a_{in}(\mathbf{X})| = O(1) \qquad \text{as } n \to \infty, \text{ a.s. } (P_\mathbf{X})$$

*and*

$$\text{(8.2)} \qquad m_n \sum_{i=1}^n a_{in}^2(\mathbf{X}) = o(1) \qquad \text{as } n \to \infty, \text{ a.s. } (P_\mathbf{X}).$$

*Then $\sum_{i=1}^n a_{in}(\mathbf{X}) g(Z(\mathbf{s}_i)) \to 0$ in $P_{\cdot|\mathbf{X}}$-probability, a.s. $(P_\mathbf{X})$.*

PROOF. This is Lemma 6.2 of [24]. □

PROOF OF THEOREM 1. Let $\xi_{n,i}^\dagger(\mathbf{t}) \equiv \int_0^1 \psi'(Z(\mathbf{s}_i) - u\mathbf{w}(\mathbf{s}_i)'(\mathbf{t} - \beta))\, du$, $1 \leq i \leq n$. By Taylor's expansion, the left-hand side of (3.3) equals

$$\text{(8.3)} \quad \begin{aligned} &\sum_{i=1}^n \mathbf{w}(\mathbf{s}_i) \psi(Z(\mathbf{s}_i)) - \sum_{i=1}^n \mathbf{w}(\mathbf{s}_i) \mathbf{w}(\mathbf{s}_i)'(\mathbf{t} - \beta) \xi_{n,i}^\dagger(\mathbf{t}) \\ &= \sum_{i=1}^n \mathbf{w}(\mathbf{s}_i) \psi(Z(\mathbf{s}_i)) - \sum_{i=1}^n \mathbf{w}(\mathbf{s}_i) \mathbf{w}(\mathbf{s}_i)'(\mathbf{t} - \beta) E\psi'(Z(\mathbf{0})) + R_n(\mathbf{t}), \end{aligned}$$

where $R_n(\mathbf{t})$ is defined by subtraction. Note that

$$\text{(8.4)} \quad \begin{aligned} \Lambda_n^{-1} R_n(\mathbf{t}) &= \Lambda_n^{-1}\left[\sum_{i=1}^n \mathbf{w}(\mathbf{s}_i) \mathbf{w}(\mathbf{s}_i)'(\mathbf{t} - \beta)\{E\psi'(Z(\mathbf{0})) - \xi_{n,i}^\dagger(\mathbf{t})\}\right] \\ &= -\left[\Lambda_n^{-1} \sum_{i=1}^n \mathbf{w}(\mathbf{s}_i) \mathbf{w}(\mathbf{s}_i)'(\Lambda_n^{-1})'\{\psi'(Z(\mathbf{s}_i)) - E\psi'(Z(\mathbf{0}))\}\right]\Lambda_n'(\mathbf{t} - \beta) \\ &\quad + \left[\Lambda_n^{-1} \sum_{i=1}^n \mathbf{w}(\mathbf{s}_i) \mathbf{w}(\mathbf{s}_i)'(\Lambda_n^{-1})'\{\psi'(Z(\mathbf{s}_i)) - \xi_{n,i}^\dagger(\mathbf{t})\}\right]\Lambda_n'(\mathbf{t} - \beta) \\ &\equiv R_{1n}(\mathbf{t}) + R_{2n}(\mathbf{t}), \qquad \text{say.} \end{aligned}$$

Next we establish convergence of $n^{-1}\sum_{i=1}^{n}\mathbf{v}(\mathbf{s}_i)\mathbf{v}(\mathbf{s}_i)'\{\psi'(Z(\mathbf{s}_i))-E\psi'(Z(\mathbf{0})))\}$ by using Lemma 1. By Hoeffding's inequality [16] and condition (C.4), it follows that for any $\varepsilon > 0$,

$$(8.5) \quad P_{\mathbf{X}}\left(\left|n^{-1}\sum_{i=1}^{n}\{\|\mathbf{v}(\mathbf{s}_i)\|^2 - E_{\mathbf{X}}\|\mathbf{v}(\mathbf{s}_i)\|^2\}\right| > \varepsilon\right) \leq C\exp(-C(\varepsilon)n/n^{2a}).$$

Hence, by the Borel–Cantelli lemma,

$$(8.6) \quad n^{-1}\sum_{i=1}^{n}(\|\mathbf{v}(\mathbf{s}_i)\|^2 - E_{\mathbf{X}}\|\mathbf{v}(\mathbf{s}_i)\|^2) = o(1) \quad \text{a.s. } (P_{\mathbf{X}}).$$

By condition (C.1),

$$(8.7) \quad E_{\mathbf{X}}\|\mathbf{v}(\mathbf{s}_1)\|^2 = \operatorname{tr}(E_{\mathbf{X}}\mathbf{v}(\mathbf{s}_1)\mathbf{v}(\mathbf{s}_1)') \to \operatorname{tr}(H) \quad \text{as } n \to \infty.$$

Hence it follows that

$$(8.8) \quad n^{-1}\sum_{i=1}^{n}\|\mathbf{v}(\mathbf{s}_i)\|^2 = O(1) \quad \text{a.s. } (P_{\mathbf{X}})$$

and

$$(8.9) \quad n^{-2}\sum_{i=1}^{n}\|\mathbf{v}(\mathbf{s}_i)\|^4 \leq \left(n^{-1}\sum_{i=1}^{n}\|\mathbf{v}(\mathbf{s}_i)\|^2\right)(n^{-1}m_{0n}^2) = o(1) \quad \text{a.s. } (P_{\mathbf{X}}).$$

Thus, by Lemma 1,

$$(8.10) \quad n^{-1}\sum_{i=1}^{n}\mathbf{v}(\mathbf{s}_i)\mathbf{v}(\mathbf{s}_i)'\{\psi'(Z(\mathbf{s}_i)) - E_{\cdot|\mathbf{X}}\psi'(Z(\mathbf{0}))\} \to 0$$

in $P_{\cdot|\mathbf{X}}$-probability, a.s. $(P_{\mathbf{X}})$.

Further, by Proposition 3.1 and Theorem 3.2 of [22] and by a conditional version of the Cramér–Wold device, it follows that for $n/\lambda_n^d \to c_1 \in (0,\infty]$,

$$(8.11) \quad L_n \equiv \eta_n^d \lambda_n^{-d/2}\sum_{i=1}^{n}\mathbf{v}(\mathbf{s}_i)\psi(Z(\mathbf{s}_i))$$

$$\xrightarrow{d} N\left(\mathbf{0}, H\sigma(\mathbf{0})/c_1 + \int \sigma(\mathbf{s})Q(\mathbf{s})\,d\mathbf{s}\right) \quad \text{a.s. } (P_{\mathbf{X}}).$$

Also, by arguments similar to those used in (8.5)–(8.7), it follows that

$$(8.12) \quad \Gamma_n \equiv n^{-1}\sum_{i=1}^{n}\mathbf{v}(\mathbf{s}_i)\mathbf{v}(\mathbf{s}_i)' \to H \quad \text{a.s. } (P_{\mathbf{X}}).$$

By (8.10)–(8.12), there exist a constant $C \in (0,\infty)$, a sequence $\varepsilon_n \downarrow 0$ and sets $B_n \in \mathcal{F}$, $n \geq 1$, such that $\lim_{n\to\infty} P_{\cdot|\mathbf{X}}(B_n) = 0$, a.s. $(P_\mathbf{X})$ and on $B_n^c$,

$$\|L_n\| \leq C\varepsilon_n^{-1/2}, \tag{8.13}$$

$$\|R_{1n}(\mathbf{t})\| \leq n\varepsilon_n \|\Lambda_n'(\mathbf{t} - \beta)\| \qquad \text{for all } \mathbf{t} \in \mathbb{R}^p \tag{8.14}$$

and

$$\|\Gamma_n - H\| \leq \varepsilon_n. \tag{8.15}$$

Without loss of generality, we assume that $\varepsilon_n^{-1} \leq \log(\lambda_n + 1)$, $n \geq 1$. Hence, on the set $B_n^c$, we may rewrite (8.3) as

$$\Lambda_n'(\mathbf{t} - \beta) = (\Gamma_n \chi_0)^{-1}[\lambda_n^{-d/2} L_n + n^{-1}\{R_{1n}(\mathbf{t}) + R_{2n}(\mathbf{t})\}].$$

By condition (C.2) and Brouwer's fixed point theorem, this implies that

$$\begin{aligned} P_{\cdot|\mathbf{X}}(\hat{\beta}_n \text{ solves } (3.3) \text{ and } \|\Lambda_n'(\hat{\beta}_n - \beta)\| \leq C\lambda_n^{-d/2}\varepsilon_n^{-1/2}) \\ = o(1) \qquad \text{a.s. } (P_\mathbf{X}). \end{aligned} \tag{8.16}$$

By (8.10) and (8.14) and condition (C.2), it follows that

$$P_{\cdot|\mathbf{X}}(\|R_n(\hat{\beta}_n)\| > \lambda_n^{-d/2}\varepsilon_n^{1/2}) = o(1) \qquad \text{a.s. } (P_\mathbf{X}). \tag{8.17}$$

Hence, by a conditional version of Slutsky's theorem (cf. Lemma 4.1 of [21]), it remains to show that

$$c_n^{-1} \sum_{i=1}^n \mathbf{w}(\mathbf{s}_i)\psi(Z(\mathbf{s}_i)) \xrightarrow{d} N(\mathbf{0}, \Sigma_c) \qquad \text{a.s. } (P_\mathbf{X}).$$

This can be established by using the Cramér–Wold device and by verifying the conditions of Theorem 3.2 and Proposition 3.1 of [22]. We omit the routine algebraic details to save space. $\square$

**9. Proof of Proposition 1.** For proving Proposition 1, we introduce some more notation. Let $S_n^{**}(\mathbf{k})$ be the sum of the DSSBB observations from the block $\mathcal{B}_n^{**}(\mathbf{k})$, $\mathbf{k} \in \mathcal{K}_n$, and let $\mathcal{S}_n^{[D]}(i)$ denote the sum of the $Y(\mathbf{s}_j)$'s for $\mathbf{s}_j \in \mathcal{B}_n^{[D]}(i)$, $i \in \mathcal{I}_n^{[D]}$. Also, let $N^{**}$ denote the DSSBB sample size, that is, $N^{**} = |\mathcal{Y}_n^{**}(R_n)|$, and let $N = E_* N^{**}$. Then

$$T_n^{**} = \sum_{\mathbf{k} \in \mathcal{K}_n} S_n^{**}(\mathbf{k})/N^{**}.$$

Since $T_n^{**}$ is a ratio estimator, we define a linearized version of it, given by

$$\tilde{T}_n^{**} = \sum_{\mathbf{k} \in \mathcal{K}_n} S_n^{**}(\mathbf{k})/N. \tag{9.1}$$

Then, by the conditional independence of the $S_n^{**}(\mathbf{k})$'s,

$$\tilde{\sigma}_n^2 \equiv \operatorname{Var}_*(\tilde{T}_n^{**}) = N^{-2}\lambda_n^d \operatorname{Var}_*\left(\sum_{\mathbf{k}\in\mathcal{K}_n} S_n^{**}(\mathbf{k})\right)$$

$$(9.2) \quad = N^{-2}\lambda_n^d|\mathcal{K}_n|\operatorname{Var}_*(S_n^{**}(\mathbf{0}))$$

$$= N^{-2}\lambda_n^d|\mathcal{K}_n|\left[|\mathcal{I}_n^{[D]}|^{-1}\sum_{i\in\mathcal{I}_n^{[D]}}(\mathcal{S}_n^{[D]}(i))^2 - \left(|\mathcal{I}_n^{[D]}|^{-1}\sum_{i\in\mathcal{I}_n^{[D]}}\mathcal{S}_n^{[D]}(i)\right)^2\right].$$

PROOF OF PROPOSITION 1. Part (i) is an immediate consequence of Lemma 5.2(ii) of [22]. This follows by noting that, in this case, the function $Q_1(\mathbf{h})$ in the lemma is given by $Q_1(\mathbf{h}) = \int f_a^2(\mathbf{s})\,d\mathbf{s}$ for all $\mathbf{h}\in\mathbb{R}^d$. To prove part (ii), note that for any $\varepsilon > 0$ and any $\delta \in (0,1)$, by (9.1),

$$P(\tilde{\sigma}_n^2 > \varepsilon) \leq P\left(N^{-2}\lambda_n^d|\mathcal{K}_n||\mathcal{I}_n^{[D]}|^{-1}\sum_{i\in\mathcal{I}_n^{[D]}}(\mathcal{S}_n^{[D]}(i))^2 > \varepsilon\right)$$

$$\leq P(|\mathcal{I}_n^{[D]}| \leq n\delta)$$

$$+ P\left(N^{-2}\lambda_n^d|\mathcal{K}_n|\sum_{i\in\mathcal{I}_n^{[D]}}(\mathcal{S}_n^{[D]}(i))^2 > \varepsilon\delta n\right)$$

$$(9.3) \quad \leq P(|\mathcal{I}_n^{[D]}| \leq n\delta)$$

$$+ \frac{N^{-2}\lambda_n^d|\mathcal{K}_n|}{\varepsilon\delta n}E\left[\sum_{i\in\mathcal{I}_n^{[D]}}(\mathcal{S}_n^{[D]}(i))^2\right]$$

$$\equiv I_1 + I_2, \quad \text{say.}$$

By using the relation $E(W) = E_\mathbf{X}(E_{\cdot|\mathbf{X}}W)$, we have

$$E\left[\sum_{i\in\mathcal{I}_n^{[D]}}(\mathcal{S}_n^{[D]}(i))^2\right]$$

$$(9.4) \quad = E\left[\sum_{i=1}^n\sum_{j=1}^n\sum_{k=1}^n Y(\mathbf{s}_j)Y(\mathbf{s}_k)\mathbb{1}(\mathbf{s}_j,\mathbf{s}_k\in\mathcal{B}_n^{[D]}(i))\mathbb{1}(\mathbf{s}_i+\mathcal{U}b_n\subset R_n)\right]$$

$$= \sum_{i=1}^n\sum_{j=1}^n\sum_{k=1}^n E_\mathbf{X}\sigma(\mathbf{s}_j-\mathbf{s}_k)\mathbb{1}(\mathbf{s}_j,\mathbf{s}_k\in\mathcal{B}_n^{[D]}(i))\mathbb{1}(\mathbf{s}_i+\mathcal{U}b_n\subset R_n).$$

Recall that $r_n = \lambda_n/b_n$. For $i \neq j \neq k$, by changing variables in the integrations below, we have

$$c_n(\mathbf{i},\mathbf{j},\mathbf{k}) \equiv E_\mathbf{X}\sigma(\mathbf{s}_j-\mathbf{s}_k)\mathbb{1}(\mathbf{s}_j,\mathbf{s}_k\in\mathcal{B}_n^{[D]}(i))\mathbb{1}(\mathbf{s}_i+\mathcal{U}b_n\subset R_n)$$

$$\leq E_{\mathbf{X}}\sigma(\lambda_n(\mathbf{X}_1 - \mathbf{X}_2))\mathbb{1}(\mathbf{X}_1, \mathbf{X}_2 \in \mathbf{X}_3 + r_n^{-1}\mathcal{U})$$

$$= \int\int\int \sigma(\lambda_n(\mathbf{x} - \mathbf{y}))\mathbb{1}(\mathbf{x}, \mathbf{y} \in \mathbf{z} + r_n^{-1}\mathcal{U}) f_a(\mathbf{x}) f_a(\mathbf{y}) f_a(\mathbf{z})\, d\mathbf{x}\, d\mathbf{y}\, d\mathbf{z}$$

$$= \lambda_n^{-d} \int\int\int \sigma(\mathbf{s})\mathbb{1}(\mathbf{y} + \lambda_n^{-1}\mathbf{s}, \mathbf{y} \in \mathbf{z} + r_n^{-1}\mathcal{U})$$

(9.5)
$$\times f_a(\mathbf{y} + \lambda_n^{-1}\mathbf{s}) f_a(\mathbf{y}) f_a(\mathbf{z})\, d\mathbf{s}\, d\mathbf{y}\, d\mathbf{z}$$

$$\leq \lambda_n^{-d} \int\int\int \sigma(\mathbf{s})\mathbb{1}(\mathbf{w} \in r_n^{-1}\mathcal{U})$$

$$\times f_a(\mathbf{w} + \mathbf{z} + \lambda_n^{-1}\mathbf{s}) f_a(\mathbf{w} + \mathbf{z}) f_a(\mathbf{z})\, d\mathbf{s}\, d\mathbf{w}\, d\mathbf{z}$$

$$= \lambda_n^{-d} r_n^{-d} \int\int\int \sigma(\mathbf{s})\mathbb{1}(\mathbf{u} \in \mathcal{U})$$

$$\times f_a(r_n^{-1}\mathbf{u} + \mathbf{z} + \lambda_n^{-1}\mathbf{s}) f_a(r_n^{-1}\mathbf{u} + \mathbf{z}) f_a(\mathbf{z})\, d\mathbf{s}\, d\mathbf{u}\, d\mathbf{z}.$$

By the continuity of $f_a(\cdot)$ and the dominated convergence theorem (DCT), the right-hand side of (9.5) is asymptotically equivalent to

(9.6) $$\lambda_n^{-d} r_n^{-d} \int \sigma(\mathbf{s})\, d\mathbf{s} \int f_a^3(\mathbf{z})\, d\mathbf{z} \equiv \Delta_n, \quad \text{say}.$$

Using similar arguments, it can be shown that

(9.7) $$N = n \int f_a^2(\mathbf{z})\, d\mathbf{z}\, (1 + o(1)) \quad \text{as } n \to \infty.$$

Hence, noting that $\mathcal{K}_n = (r_n + 1)^d$, by (4.6), (9.3)–(9.7), it follows that

$$I_2 \leq \frac{N^{-2}\lambda_n^d(r_n + 1)^d}{\varepsilon\delta n}\left[cn^2\sigma(\mathbf{0}) + \sum_{\mathbf{i}\neq\mathbf{j}\neq\mathbf{k}} c_n(\mathbf{i}, \mathbf{j}, \mathbf{k})\right]$$

(9.8)
$$\leq C\frac{\lambda_n^{2d}}{nb_n^d} + \frac{N^{-2}\lambda_n^d(r_n + 1)^d}{\varepsilon\delta n}n^3\Delta_n(1 + o(1))$$

$$= (\varepsilon\delta)^{-1}\left(\int \sigma(\mathbf{s})\, d\mathbf{s}\right)\left(\int f_a^3(\mathbf{z})\, d\mathbf{z}\right)(n^2 N^{-2})(1 + o(1))$$

$$= (\varepsilon\delta)^{-1}\sigma_\infty^2\left[\int f_a^3(\mathbf{z})\, d\mathbf{z} \bigg/ \left(\int f_a^2(\mathbf{z})\, d\mathbf{z}\right)^3\right](1 + o(1)).$$

Note that $\int f_a^3(\mathbf{z})\, d\mathbf{z} = \int g_a(\mathbf{x}_1)^3\, d\mathbf{x}_1 < 2[(\frac{a}{4})^3\frac{2}{a} + b^3(\frac{1}{2} - \frac{2}{a})]$ while $(\int f_a^2(\mathbf{z})\, d\mathbf{z})^3 = [\int g_a(\mathbf{x}_1)^2\, d\mathbf{x}_1]^3 > [2(\frac{a}{4})^2\frac{1}{a}]^3 = \frac{a^3}{8^3} \gg \int f_a^3(\mathbf{z})\, d\mathbf{z}$, as $a \to \infty$, where $b = \frac{a}{4(a-3)} < \frac{a}{4}$ for all $a > 4$. Hence there exist $a > 4$ and $\delta = \delta(a) \in (0,1)$, $\eta(a) \in (0,1)$ (with both $\delta$ and $\eta$ sufficiently close to 1) such that with $\varepsilon = \eta\sigma_\infty^2$, the leading term on the right-hand side of (9.8) is $< 1$.

To obtain an estimate of $I_1$, define $W_{in} = \mathbb{1}(\lambda_n \mathbf{X}_i + b_n \mathcal{U} \subset R_n)$, $1 \leq i \leq n$. Then, by the DCT, $EW_{1n} \to 1$ as $n \to \infty$. Hence, for any $0 < \delta < 1$, by Chebyshev's inequality,

$$\lim_{n\to\infty} I_1 = \lim_{n\to\infty} P(|\mathcal{I}_n^{[D]}| \leq n\delta) = \lim_{n\to\infty} P\left(\sum_{i=1}^n W_{in} \leq n\delta\right)$$
$$\leq \lim_{n\to\infty} P\left(\left|\sum_{i=1}^n (W_{in} - EW_{in})\right| > n|EW_{1n} - \delta|\right)$$
$$\leq \lim_{n\to\infty} \operatorname{Var}\left(\sum_{i=1}^n W_{in}\right) \bigg/ [n^2(EW_{1n} - \delta)^2] = 0.$$

Thus, it follows that (4.12) holds with $\hat{\sigma}_n^2$ replaced by $\tilde{\sigma}_n^2$. To complete the proof, it is enough to show that

$$\hat{\sigma}_n^2 - \tilde{\sigma}_n^2 = o_p(1).$$

This can be shown by using the identity $T_n^{**} = \tilde{T}_n^{**} + (N - N^{**})\tilde{T}_n^{**}/N^{**}$, the Cauchy–Schwarz inequality and the boundedness of the variables $Y(\cdot)$'s. We omit the details. $\square$

## 10. Proof of Theorem 2.

PROOF OF THEOREM 2. For a positive definite matrix $\Sigma$, let $\Phi(\cdot; \Sigma)$ denote the probability measure corresponding to the $N(\mathbf{0}, \Sigma)$ distribution. By Taylor's expansion, from (5.11), for any $\mathbf{t} \in \mathbb{R}^p$, we get

$$0 = \sum_{\mathbf{k} \in \mathcal{K}_n} [S_n^*(\mathbf{k}; \mathbf{t}) - \hat{c}_n(\mathbf{k})]$$
$$= \sum_{\mathbf{k} \in \mathcal{K}_n} [S_n^*(\mathbf{k}; \hat{\beta}_n) - \hat{c}_n(\mathbf{k})]$$
$$(10.1) \quad + \sum_{\mathbf{k} \in \mathcal{K}_n} \sum_{i=1}^n \bigg[ \mathbf{w}(\mathbf{s}_i)\mathbf{w}(\mathbf{s}_i)'(\mathbf{t} - \hat{\beta}_n)\mathbb{1}(\mathbf{s}_i \in B_n(I_{\mathbf{k}}; \mathbf{k}))$$
$$\times \int_0^1 \psi(\hat{Z}(\mathbf{s}_i) - u\mathbf{w}(\mathbf{s}_i)'(\mathbf{t} - \hat{\beta}_n))\,du \bigg]$$
$$\equiv \sum_{\mathbf{k} \in \mathcal{K}_n} [S_n^*(\mathbf{k}; \hat{\beta}_n) - \hat{c}_n(\mathbf{k})] + \Lambda_n \Gamma_n \Lambda_n'(\mathbf{t} - \hat{\beta}_n)n\chi_0 + R_n^*(\mathbf{t}), \quad \text{say},$$

where the remainder term is defined by subtraction.

Now, by arguments similar to those employed in the proof of Theorem 1, it is enough to prove that there exists a nonrandom sequence $\{\varepsilon_n\}_{n\geq 1} \subset (0,1)$ with $\varepsilon_n \downarrow 0$ as $n \to \infty$ such that for any $\varepsilon_0 \in (0,1)$,

$$
\text{(10.2)} \quad P_{\cdot|\mathbf{X}}\left(\sup_{B \in \mathcal{C}}\left|P_*\left(C_n^{-1}\sum_{\mathbf{k} \in \mathcal{K}_n}[S_n^*(\mathbf{k};\hat{\beta}_n) - \hat{c}_n(\mathbf{k})] \in B\right) - \Phi(B;\Sigma_c)\right| > \varepsilon_0\right)
$$
$$
= o(1), \quad n \to \infty, \text{ a.s. } (P_\mathbf{X}),
$$

and

$$
\text{(10.3)} \quad \begin{aligned} P_{\cdot|\mathbf{X}}(P_*(\|C_n^{-1}R_n^*(\mathbf{t})\| > \varepsilon_n \|V_n(\mathbf{t})\| \\ \text{for some } \mathbf{t} \in \mathbb{R}^p \text{ with } \|V_n(\mathbf{t})\| \leq |\log \varepsilon_n|) > \varepsilon_0) \\ = o(1), \quad n \to \infty, \text{ a.s. } (P_\mathbf{X}), \end{aligned}
$$

where we recall that $C_n^{-1} = \lambda_n^{-d/2}\eta_n^d\Lambda_n^{-1}$, and where $V_n(\mathbf{t}) = \lambda_n^{d/2}\Lambda_n'(\mathbf{t} - \hat{\beta}_n)$.

Consider (10.2) first. Let $\hat{\Sigma}_n = \sum_{\mathbf{k} \in \mathcal{K}_n} \text{Var}_*(C_n^{-1}S_n^*(\mathbf{k};\hat{\beta}_n))$ and let $\hat{\delta}_n$ be the smallest eigenvalue of $\hat{\Sigma}_n$. By Theorem 17.2 of [3], there exists a constant $C(p) \in (0,\infty)$ such that, on the set $\{\hat{\delta}_n > 0\}$,

$$
\text{(10.4)} \quad \begin{aligned} \sup_{B \in \mathcal{C}}\left|P_*\left(C_n^{-1}\sum_{\mathbf{k} \in \mathcal{K}_n}[S_n^*(\mathbf{k};\hat{\beta}_n) - \hat{c}_n(\mathbf{k})] \in B\right) - \Phi(B;\hat{\Sigma}_n)\right| \\ \leq C(p)\sum_{\mathbf{k} \in \mathcal{K}_n}E_*\|C_n^{-1}S_n^*(\mathbf{k};\hat{\beta}_n)\|^3\hat{\delta}_n^{-3}. \end{aligned}
$$

In Lemma 2, we show that

$$
\text{(10.5)} \quad \|\hat{\Sigma}_n - \Sigma_c\| \to 0 \quad \text{in } P_{\cdot|\mathbf{X}}\text{- probability, a.s. } (P_\mathbf{X}).
$$

Also, from part (b) of Lemma 2, it follows that

$$
\text{(10.6)} \quad \sum_{\mathbf{k} \in \mathcal{K}_n}E_*\|C_n^{-1}S_n^*(\mathbf{k};\hat{\beta}_n)\|^3 \to 0 \quad \text{in } P_{\cdot|\mathbf{X}}\text{- probability, a.s. } (P_\mathbf{X}).
$$

Hence, by (10.4)–(10.6), (10.2) follows.

Next consider (10.3). Note that $n^{-1}\sum_{i=1}^n \mathbf{w}(\mathbf{s}_i)\mathbf{w}(\mathbf{s}_i)' = \Lambda_n\Gamma_n\Lambda_n'$. Hence, by condition (C.2), for any $\mathbf{t} \in \mathbb{R}^p$, we have

$\|C_n^{-1}R_n^*(\mathbf{t})\|$

$$
\leq \left\|n^{-1}\sum_{\mathbf{k} \in \mathcal{K}_n}\sum_{i=1}^n \mathbf{v}(\mathbf{s}_i)\mathbf{v}(\mathbf{s}_i)'\mathbb{1}(\mathbf{s}_i \in B_n(I_\mathbf{k};\mathbf{k})) - \Gamma_n\right\|\|V_n(\mathbf{t})\|\|\chi_0\|
$$
$$
+ \frac{C}{n}\sum_{\mathbf{k} \in \mathcal{K}_n}\sum_{i=1}^n \|\mathbf{v}(\mathbf{s}_i)\|^{2+\gamma}\mathbb{1}(\mathbf{s}_i \in B_n(I_\mathbf{k};\mathbf{k}))
$$
$$
\text{(10.7)} \quad \times (\|V_n(\mathbf{t})\|^\gamma + \|V_n(\beta)\|^\gamma)\lambda_n^{-d\gamma/2}\|V_n(\mathbf{t})\|
$$

$$+ \left\| n^{-1} \sum_{\mathbf{k} \in \mathcal{K}_n} \sum_{i=1}^{n} \mathbf{v}(\mathbf{s}_i)\mathbf{v}(\mathbf{s}_i)'\{\psi'(Z(\mathbf{s}_i)) - \chi_0\}\mathbb{1}(\mathbf{s}_i \in B_n(I_\mathbf{k}; \mathbf{k})) \right\| \|V_n(\mathbf{t})\|$$

$$\equiv R_{1n}^*(\mathbf{t}) + R_{2n}^*(\mathbf{t}) + R_{3n}^*(\mathbf{t}), \quad \text{say.}$$

Consider $R_{1n}^*(\mathbf{t})$ first. For any $\varepsilon > 0$ and $A \subset \mathbb{R}^d$, let $A^{-\varepsilon} = \{x \in A : x + [-1,1]^d \varepsilon \subset A\}$. Let $u_n = b_n/\lambda_n$, $n \geq 1$. Then it is easy to show that uniformly in $x \in R_0^{-u_0}$,

$$(10.8) \qquad |\{\mathbf{m} \in \mathcal{I}_n : \lambda_n \mathbf{x} \in \mathbf{m} + [0,1)^d\}| = b_n^d + O(b_n^{d-1}).$$

Hence, by (10.8) and the fact that $|\mathcal{K}_{1n}| \sim \text{vol}(R_0)\lambda_n^d b_n^{-d}$ and $|\mathcal{I}_n| \sim \text{vol}(R_0)\lambda_n^d$, we have

$$\tilde{\Gamma}_n \equiv E_* n^{-1} \sum_{\mathbf{k} \in \mathcal{K}_{1n}} \sum_{j=1}^{n} \mathbf{v}(\mathbf{s}_j)\mathbf{v}(\mathbf{s}_j)'\mathbb{1}(\mathbf{s}_j \in B_n(I_\mathbf{k}; \mathbf{k}))\mathbb{1}(\mathbf{s}_j \in \lambda_n R_0^{-u_n})$$

$$(10.9) \quad = n^{-1} \sum_{j=1}^{n} \mathbf{v}(\mathbf{s}_j)\mathbf{v}(\mathbf{s}_j)'\mathbb{1}(\mathbf{s}_j \in \lambda_n R_0^{-u_n})|\mathcal{K}_{1n}||\mathcal{I}_n|^{-1} \sum_{\mathbf{i} \in \mathcal{I}_n} \mathbb{1}(\mathbf{s}_j \in B_n(\mathbf{i}; \mathbf{0}))$$

$$= \left[ n^{-1} \sum_{j=1}^{n} \mathbf{v}(\mathbf{s}_j)\mathbf{v}(\mathbf{s}_j)'\mathbb{1}(\mathbf{s}_j \in \lambda_n R_0^{-u_n}) \right](1 + o(1)).$$

Further, note that by the boundary condition on $R_0$ and (8.6),

$$(10.10) \quad E_*\left[ n^{-1} \sum_{j=1}^{n} \|\mathbf{v}(\mathbf{s}_j)\|^2 \sum_{\mathbf{k} \in \mathcal{K}_{2n}} \mathbb{1}(\mathbf{s}_j \in B_n(I_\mathbf{k}; \mathbf{k})) \right]$$

$$\leq n^{-1} \sum_{j=1}^{n} \|\mathbf{v}(\mathbf{s}_j)\|^2 |\mathcal{K}_{2n}||\mathcal{I}_n|^{-1}b_n^d = O(u_n) \quad \text{a.s.} \ (P_\mathbf{X}).$$

Next, note that by (8.8) and the condition on $m_{0n}$, for any $a \in (0, \infty)$,

$$E_\mathbf{X}\|\mathbf{v}(\mathbf{s}_1)\|^2 \mathbb{1}(\mathbf{s}_1 \notin \lambda_n R_0^{-u_n})$$

$$(10.11) \quad \leq (E_\mathbf{X}\|\mathbf{v}(\mathbf{s}_1)\|^{2+a})^{1/(2+a)}[P_\mathbf{X}(\mathbf{X}_1 \notin R_0^{-u_n})]^{a/(2+a)}$$

$$\leq m_{0n}^{a/(2+a)}(E_\mathbf{X}\|\mathbf{v}(\mathbf{s}_1)\|^2)^{a/(2+a)} \cdot O(u_n^{a/(2+a)})$$

$$= o(1) \quad \text{as } n \to \infty.$$

Hence by arguments leading to (8.6), a.s. $(P_\mathbf{X})$,

$$(10.12) \quad E_*\left\{ n^{-1} \sum_{j=1}^{n} \|\mathbf{v}(\mathbf{s}_j)\|^2 \mathbb{1}(\mathbf{s}_j \in \lambda_n R_0^{-u_n}) \sum_{\mathbf{k} \in \mathcal{K}_{1n}} \mathbb{1}(\mathbf{s}_j \in B_n(I_\mathbf{k}; \mathbf{k})) \right\}$$

$$\leq n^{-1} \sum_{j=1}^{n} \|\mathbf{v}(\mathbf{s}_j)\|^2 \mathbb{1}(\mathbf{s}_j \in \lambda_n R_0^{-u_n})|\mathcal{K}_{1n}||\mathcal{I}_n|^{-1}b_n^d = o(1).$$

Now, by (10.9), (10.10), (10.12) and Chebyshev's inequality, there exists $\varepsilon_{1n} \downarrow 0$ such that (10.3) holds with $\varepsilon_n = \varepsilon_{1n}$ and $R_n^*(\mathbf{t})$ replaced by $R_{1n}^*(\mathbf{t})$. Also, in view of condition (C.2), (10.3) holds for $R_{2n}^*(\mathbf{t})$ by Markov's inequality [of the form $P_*(|W^*| > \varepsilon_{2n}) \leq E_*|W^*|/\varepsilon_{2n}$] with $\varepsilon_n = \varepsilon_{2n}$, where $\varepsilon_{2n} = \lambda_n^{-d/32}$, say. Hence, it remains to prove (10.3) for $R_{3n}^*(\mathbf{t})$. To that end, let $W_{jl}(\mathbf{s}_i) = \mathbf{e}_j' \mathbf{v}(\mathbf{s}_i)' \mathbf{v}(\mathbf{s}_i) \mathbf{e}_l [\psi'(Z(\mathbf{s}_i)) - \chi_0]$, $1 \leq j, l \leq p$, $i = 1, \ldots, n$. Then, by the independence of the resampled blocks,

$$
\begin{aligned}
E_{\cdot|\mathbf{X}} &\left[ \mathrm{Var}_* \left( \sum_{\mathbf{k} \in \mathcal{K}_n} n^{-1} \sum_{i=1}^n W_{jl}(\mathbf{s}_i) \mathbb{1}(\mathbf{s}_i \in B_n(I_\mathbf{k}; \mathbf{k})) \right) \right] \\
&\leq E_{\cdot|\mathbf{X}} \left[ \sum_{\mathbf{k} \in \mathcal{K}_n} E_* \left\{ n^{-1} \sum_{i=1}^n W_{jl}(\mathbf{s}_i) \mathbb{1}(\mathbf{s}_i \in B_n(I_\mathbf{k}; \mathbf{k})) \right\}^2 \right] \\
&\leq E_{\cdot|\mathbf{X}} \left[ C n^{-2} |\mathcal{I}_n|^{-1} \sum_{\mathbf{k} \in \mathcal{K}_n} \sum_{\mathbf{i} \in \mathcal{I}_n} \left| \sum_{i=1}^n W_{jl}(\mathbf{s}_i) \mathbb{1}(\mathbf{s}_i \in B_n(I_\mathbf{k}; \mathbf{k})) \right|^2 \right] \\
&= O(n^{-2} b_n^{2d} (m_{0n}^2 m_n)^2) = o(1) \qquad \text{a.s. } (P_\mathbf{X})
\end{aligned}
$$
(10.13)

by Lemma 2(b). Hence (10.3) holds for $R_{3n}^*(\mathbf{t})$. This completes the proof of Theorem 2. □

LEMMA 2. *Assume that $b_n$ satisfies (5.15). Let $g:\mathbb{R} \to \mathbb{R}$ be a Borel-measurable function satisfying $Eg(Z(\mathbf{0})) = 0$ and let $\{w_{ni}(\mathbf{X}) : 1 \leq i \leq n\}_{n \geq 1}$ be $\chi$-measurable weight variables with $M_n(\mathbf{X}) \equiv \max\{|w_{ni}(\mathbf{X})| : i = 1, \ldots, n\} < \infty$. Suppose that $E|g(Z(\mathbf{0}))|^{2r+\delta} < \infty$ and $A(r, \delta) < \infty$ for some $r \in \mathbb{N}$ and $\delta \in (0, \infty)$. Then*

$$
\sum_{\mathbf{k} \in \mathcal{K}_n} \sum_{\mathbf{i} \in \mathcal{I}_n} E_{\cdot|\mathbf{X}} \left| \sum_{j=1}^n w_{nj}(\mathbf{X}) g(Z(\mathbf{s}_j)) \mathbb{1}(\mathbf{s}_j \in B_n(\mathbf{i}; \mathbf{k})) \right|^{2r}
$$
(10.14)
$$
= O(|\mathcal{K}_n| |\mathcal{I}_n| b_n^{dr} [M_n(\mathbf{X}) m_n]^{2r}) \qquad a.s. \ (P_\mathbf{X}),
$$

*and for any $B_n \subset R_n$,*

$$
E_{\cdot|\mathbf{X}} \left( \sum_{i=1}^n w_{ni}(\mathbf{X}) g(Z(\mathbf{s}_j)) \mathbb{1}(\mathbf{s}_j \in B_n) \right)^{2r}
$$
(10.15)
$$
= O(|J_n|^r [M_n(\mathbf{X}) m_n]^{2r}) \qquad a.s. \ (P_\mathbf{X}),
$$

*where $J_n = \{\mathbf{j} \in \mathbb{Z}^d : B_n \cap \{\mathbf{j} + [0, 1)^d\} \neq \varnothing\}$.*

PROOF. Using the exponential inequality in Lemma 5.1 of [22], it can be shown that there exists a constant $C \in (0, \infty)$ such that

$$P_{\mathbf{X}}\left(\max_{\|\mathbf{j}\|_1 \leq d\lambda_n/2, \mathbf{i} \in \mathcal{I}_n, \mathbf{k} \in \mathcal{K}_n} \sum_{i=1}^n \mathbb{1}(\lambda_n \mathbf{X}_i \in \{\mathbf{j} + [0,1)^d\} \cap B_n(\mathbf{i}; \mathbf{k})) > Cm_n \text{ i.o.}\right)$$
(10.16)
$$= 0.$$

Next, grouping the terms in the sum $\sum_{j=1}^n w_{nj}(\mathbf{X})g(Z(\mathbf{s}_j))\mathbb{1}(\mathbf{s}_j \in B_n(\mathbf{i};\mathbf{k}))$ corresponding to the $\mathbf{s}_i$'s in the unit cubes $\{\mathbf{j} + [0,1)^d\} \cap B_n(\mathbf{i};\mathbf{k})$'s and applying Theorem 1.4.1.1 of [9] for each $\mathbf{i} \in \mathcal{I}_n, \mathbf{k} \in \mathcal{K}_n$, one gets (10.14). The proof of (10.15) is similar and is omitted. $\square$

LEMMA 3. *Suppose the conditions of Theorem 2 hold. Then*
$$\|\hat{\Sigma}_n - \Sigma_c\| \to 0 \quad \text{in } P_{\cdot|\mathbf{X}}\text{- probability, a.s. } (P_{\mathbf{X}}).$$

PROOF. For notational simplicity, we prove Lemma 3 assuming that $p = 1$. The proof for general $p \geq 1$ is similar, except for some additional complexity in notation. Let $\tilde{S}(\mathbf{i},\mathbf{k}) = \sum_{j=1}^n \mathbf{w}(\mathbf{s}_j)\psi(Z(\mathbf{s}_j))\mathbb{1}(\mathbf{s}_j \in B_n(\mathbf{i};\mathbf{k})), \mathbf{i} \in \mathcal{I}_n, \mathbf{k} \in \mathcal{K}_n$, and $\tilde{\Sigma}_n = \sum_{\mathbf{k} \in \mathcal{K}_n}[|\mathcal{I}_n|^{-1}\sum_{\mathbf{i} \in \mathcal{I}_n}(C_n^{-1}\tilde{S}(\mathbf{i},\mathbf{k}))^2 - (|\mathcal{I}_n|^{-1}\sum_{\mathbf{i} \in \mathcal{I}_n}(C_n^{-1}\tilde{S}(\mathbf{i},\mathbf{k})))^2]$. Then, by condition (C.2), Lemma 1 and (8.16) and (10.16), it follows that

(10.17) $\quad \hat{\Sigma}_n - \tilde{\Sigma}_n \to 0 \quad \text{in } P_{\cdot|\mathbf{X}}\text{-probability, a.s. } (P_{\mathbf{X}}).$

Next, using the regrouping arguments on pages 298–299 of [21], and by (10.14) of Lemma 2, it follows that

$$\sum_{\mathbf{k} \in \mathcal{K}_n} E_{\cdot|\mathbf{X}}\left(\sum_{\mathbf{i} \in \mathcal{I}_n}[\{\Lambda_n^{-1}\tilde{S}(\mathbf{i},\mathbf{k})\}^2 - \{E_{\cdot|\mathbf{X}}(\Lambda_n^{-1}\tilde{S}(\mathbf{i},\mathbf{k}))\}^2]\right)^2$$
(10.18)
$$= O(|\mathcal{K}_n|\lambda_n^d b_n^d m_{0n}^4 m_n^4 b_n^{2d}) \quad \text{a.s. } (P_{\mathbf{X}}),$$

and by (10.15) of Lemma 2,

(10.19) $\quad \sum_{\mathbf{k} \in \mathcal{K}_n} E_{\cdot|\mathbf{X}}\left(|\mathcal{I}_n|^{-1}\sum_{\mathbf{i} \in \mathcal{I}_n}\Lambda_n^{-1}\tilde{S}(\mathbf{i},\mathbf{k})\right)^2 = O(b_n^d m_{0n}^2 m_n^2) \quad \text{a.s. } (P_{\mathbf{X}}).$

By (10.17)–(10.19), it now remains to show that

(10.20) $\quad \lim_{n \to \infty} \sum_{\mathbf{k} \in \mathcal{K}_n} |\mathcal{I}_n|^{-1} \sum_{\mathbf{i} \in \mathcal{I}_n} E_{\cdot|\mathbf{X}}(C_n^{-1}\tilde{S}(\mathbf{i},\mathbf{k}))^2 = \Sigma_c \quad \text{a.s. } (P_{\mathbf{X}}).$

Let $\tilde{\Sigma}_{jn} = \sum_{\mathbf{k} \in \mathcal{K}_{jn}} |\mathcal{I}_n|^{-1} \sum_{\mathbf{i} \in \mathcal{I}_n} E_{\cdot|\mathbf{X}}(C_n^{-1}\tilde{S}(\mathbf{i},\mathbf{k}))^2, j = 1, 2$. Then

$$E_{\mathbf{X}}(\tilde{\Sigma}_{1n}) = n^{-2}\lambda_n^d |\mathcal{K}_{1n}||\mathcal{I}_n|^{-1}$$
$$\times \left[\sum_{\mathbf{i} \in \mathcal{I}_n} \sigma_\psi(\mathbf{0})[nE_{\mathbf{X}}v(\lambda_n \mathbf{X}_1)^2\right.$$

$$+ \sum_{\mathbf{i} \in \mathcal{I}_n} \{n(n-1)E_{\mathbf{X}} v(\lambda_n \mathbf{X}_1) v(\lambda_n \mathbf{X}_2)$$

$$\times \sigma_\psi(\lambda_n(\mathbf{X}_1 - \mathbf{X}_2))\mathbb{1}(\lambda_n \mathbf{X}_1, \lambda_n \mathbf{X}_2 \in B_n(\mathbf{i}; \mathbf{0}))\}\Bigg]$$

$$\equiv \Sigma_{11n} + \Sigma_{12n}, \qquad \text{say,}$$

where $\sigma_\psi(\mathbf{s}) = E\psi(Z(\mathbf{0}))\psi(Z(\mathbf{s}))$, $\mathbf{s} \in \mathbb{R}^d$. Note that by (10.8) and condition (C.1) [cf. (3.6)],

$$\Sigma_{11n} = \frac{\lambda_n^d |\mathcal{K}_{1n}|}{n|\mathcal{I}_n|} \sigma_\psi(\mathbf{0}) \Bigg[ \int_{R_0^{-u_n}} v(\lambda_n \mathbf{x})^2 f(\mathbf{x}) \sum_{\mathbf{i} \in \mathcal{I}_n} \mathbb{1}(\lambda_n \mathbf{x}_1 \in \mathbf{i} + [0,1)^d b_n)\,d\mathbf{x}$$

(10.21)
$$+ O\Bigg(\int_{R_0 \setminus R_0^{-u_n}} v(\lambda_n \mathbf{x})^2 f(\mathbf{x})\,d\mathbf{x}\, b_n^d\Bigg) \Bigg]$$

$$= \sigma_\psi(\mathbf{0}) H n^{-1} \lambda_n^d (1 + o(1)) \qquad \text{as } n \to \infty.$$

It can be shown that uniformly in $\mathbf{x} \in R_o^{-u_n}$ and $\|\mathbf{s}\| \leq b_n^{1/2}$, $|\{\mathbf{i} \in \mathcal{I}_n : \lambda_n \mathbf{x} \in B_n(\mathbf{i}; \mathbf{0}), \lambda_n \mathbf{x} + \mathbf{s} \in B_n(\mathbf{i}; \mathbf{0})\}| = b_n^d(1 + O(b_n^{-1/2}))$ as $n \to \infty$. Hence, using condition (C.1) [cf. (3.7)], condition (C.5), and change of variables with $\tilde{R}_0 \equiv \{\mathbf{x} - \mathbf{y} : \mathbf{x}, \mathbf{y} \in R_0\}$, we get

$$\Sigma_{12n} = n^{-1} \lambda_n^d |\mathcal{K}_{1n}||\mathcal{I}_n|^{-1}(n-1)$$

$$\times \int_{\tilde{R}_0} \int_{R_0 \cap R_0 + \mathbf{u}} \sigma_\psi(\lambda_n \mathbf{u}) v(\lambda_n \mathbf{x}) v(\lambda_n (\mathbf{x} - \mathbf{u}))$$

$$\times \Bigg\{\sum_{\mathbf{i} \in \mathcal{I}_n} \mathbb{1}(\lambda_n \mathbf{x}, \lambda_n(\mathbf{x} - \mathbf{u}) \in B_n(\mathbf{i}; \mathbf{0}))\Bigg\} f(\mathbf{x}) f(\mathbf{x} - \mathbf{u})\,d\mathbf{x}\,d\mathbf{u}$$

(10.22)
$$= \frac{|\mathcal{K}_{1n}|}{|\mathcal{I}_n|} \int_{\|\mathbf{s}\|^2 \leq b_n} \sigma_\psi(\mathbf{s}) \int_{R_0^{-u_n} \cap (R_0 + \lambda_n^{-1} \mathbf{s})} v(\lambda_n \mathbf{x}) v(\lambda_n \mathbf{x} - \mathbf{s})$$

$$\times f(\mathbf{x}) f(\mathbf{x} - \lambda_n^{-1} \mathbf{s})$$

$$\times \Bigg\{\sum_{\mathbf{i} \in \mathcal{I}_n} \mathbb{1}(\lambda_n \mathbf{x}, \lambda_n \mathbf{x} - \mathbf{s} \in B_n(\mathbf{i}; \mathbf{0}))\Bigg\} d\mathbf{x}\,d\mathbf{s}$$

$$+ o(1)$$

$$= \int \sigma_\psi(\mathbf{s}) Q(\mathbf{s})\,d\mathbf{s} + o(1) \qquad \text{as } n \to \infty.$$

Thus, by (10.21)–(10.22), $E_{\mathbf{X}} \tilde{\Sigma}_{1n} \to \Sigma_c$ as $n \to \infty$. Since $|\mathcal{K}_{2n}| = o(|\mathcal{K}_{1n}|)$ and $B_n(\mathbf{i}; \mathbf{k}) \subset B_n(\mathbf{i}; \mathbf{0})$ for all $\mathbf{k} \in \mathcal{K}_{2n}$, a similar set of arguments shows that

$E_{\mathbf{X}}\tilde{\Sigma}_{2n} = o(1)$ as $n \to \infty$. To complete the proof of (10.20), it is now enough to show that

$$\tilde{\Delta}_n \equiv \tilde{\Sigma}_{1n} + \tilde{\Sigma}_{2n} - (E_{\mathbf{X}}\tilde{\Sigma}_{1n} + E_{\mathbf{X}}\tilde{\Sigma}_{2n}) \to 0 \qquad \text{a.s. } (P_{\mathbf{X}}).$$

This can be done by writing $\tilde{\Delta}_n$ as a U-statistic of order 2 in $\{\mathbf{X}_1, \ldots, \mathbf{X}_n\}$ and using the arguments developed in the proof of Lemma 5.2 of [22]. We omit the routine details. $\square$

## REFERENCES

[1] ANDERSON, D. A., TURNER, M. G., FORESTER, J. D., ZHU, J., BOYCE, M. S., BEYER, H. and STOWELL, L. (2005). Scale-dependent summer resource selection by reintroduced elk in Wisconsin, USA. *J. Wildlife Management* **69** 298–310.

[2] ANDERSON, T. W. (1971). *The Statistical Analysis of Time Series*. Wiley, New York. MR0283939

[3] BHATTACHARYA, R. N. and RANGA RAO, R. (1986). *Normal Approximation and Asymptotic Expansions*. Krieger, Melbourne, FL. MR0855460

[4] BILLINGSLEY, P. (1968). *Convergence of Probability Measures*. Wiley, New York. MR0233396

[5] BRADLEY, R. C. (1989). A caution on mixing conditions for random fields. *Statist. Probab. Lett.* **8** 489–491. MR1040812

[6] BRADLEY, R. C. (1993). Equivalent mixing conditions for random fields. *Ann. Probab.* **21** 1921–1926. MR1245294

[7] CARLSTEIN, E. (1986). The use of subseries values for estimating the variance of a general statistic from a stationary sequence. *Ann. Statist.* **14** 1171–1179. MR0856813

[8] CRESSIE, N. (1993). *Statistics for Spatial Data*, rev. ed. Wiley, New York. MR1239641

[9] DOUKHAN, P. (1994). *Mixing. Properties and Examples. Lecture Notes in Statist.* **85**. Springer, New York. MR1312160

[10] FREEDMAN, D. A. (1981). Bootstrapping regression models. *Ann. Statist.* **9** 1218–1228. MR0630104

[11] GARCIA-SOIDAN, P. H. and HALL, P. (1997). On sample reuse methods for spatial data. *Biometrics* **53** 273–281. MR1450185

[12] GRENANDER, U. (1954). On estimation of regression coefficients in the case of an autocorrelated disturbance. *Ann. Math. Statist.* **25** 252–272. MR0062402

[13] HALL, P. (1985). Resampling a coverage pattern. *Stochastic Process. Appl.* **20** 231–246. MR0808159

[14] HALL, P., HOROWITZ, J. L. and JING, B.-Y. (1995). On blocking rules for the bootstrap with dependent data. *Biometrika* **82** 561–574. MR1366282

[15] HALL, P. and PATIL, P. (1994). Properties of nonparametric estimators of autocovariance for stationary random fields. *Probab. Theory Related Fields* **99** 399–424. MR1283119

[16] HOEFFDING, W. (1963). Probability inequalities for sums of bounded random variables. *J. Amer. Statist. Assoc.* **58** 13–30. MR0144363

[17] KÜNSCH, H. R. (1989). The jackknife and the bootstrap for general stationary observations. *Ann. Statist.* **17** 1217–1241. MR1015147

[18] LAHIRI, S. N. (1992). Bootstrapping $M$-estimators of a multiple linear regression parameter. *Ann. Statist.* **20** 1548–1570. MR1186265

[19] LAHIRI, S. N. (1993). On the moving block bootstrap under long range dependence. *Statist. Probab. Lett.* **18** 405–413. MR1247453

[20] LAHIRI, S. N. (1999). Asymptotic distribution of the empirical spatial cumulative distribution function predictor and prediction bands based on a subsampling method. *Probab. Theory Related Fields* **114** 55–84. MR1697139

[21] LAHIRI, S. N. (2003). *Resampling Methods for Dependent Data*. Springer, New York. MR2001447

[22] LAHIRI, S. N. (2003). Central limit theorems for weighted sums of a spatial process under a class of stochastic and fixed designs. *Sankhyā Ser. A* **65** 356–388. MR2028905

[23] LAHIRI, S. N., KAISER, M. S., CRESSIE, N. and HSU, N.-J. (1999). Prediction of spatial cumulative distribution functions using subsampling (with discussion). *J. Amer. Statist. Assoc.* **94** 86–110. MR1689216

[24] LAHIRI, S. N. and MUKHERJEE, K. (2004). Asymptotic distributions of $M$-estimators in a spatial regression model under some fixed and stochastic spatial sampling designs. *Ann. Inst. Statist. Math.* **56** 225–250. MR2067154

[25] LEE, Y.-D. and LAHIRI, S. N. (2002). Least squares variogram fitting by spatial subsampling. *J. R. Stat. Soc. Ser. B Stat. Methodol.* **64** 837–854. MR1979390

[26] LIU, R. Y. and SINGH, K. (1992). Moving blocks jackknife and bootstrap capture weak dependence. In *Exploring the Limits of Bootstrap* (R. LePage and L. Billard, eds.) 225–248. Wiley, New York. MR1197787

[27] NORDMAN, D. and LAHIRI, S. N. (2004). On optimal spatial subsample size for variance estimation. *Ann. Statist.* **32** 1981–2027. MR2102500

[28] PARTHASARATHY, K. (1967). *Probability Measures on Metric Spaces*. Academic Press, New York. MR0226684

[29] POLITIS, D. N., PAPARODITIS, E. and ROMANO, J. P. (1998). Large sample inference for irregularly spaced dependent observations based on subsampling. *Sankhyā Ser. A* **60** 274–292. MR1711685

[30] POLITIS, D. N., PAPARODITIS, E. and ROMANO, J. P. (1999). Resampling marked point processes. In *Multivariate Analysis, Design of Experiments, and Survey Sampling: A Tribute to J. N. Srivastava* (S. Ghosh, ed.) 163–185. Dekker, New York. MR1719074

[31] POLITIS, D. N. and ROMANO, J. P. (1993). Nonparametric resampling for homogeneous strong mixing random fields. *J. Multivariate Anal.* **47** 301–328. MR1247380

[32] POLITIS, D. N. and ROMANO, J. P. (1994). The stationary bootstrap. *J. Amer. Statist. Assoc.* **89** 1303–1313. MR1310224

[33] POLITIS, D. N. and ROMANO, J. P. (1994). Large sample confidence regions based on subsamples under minimal assumptions. *Ann. Statist.* **22** 2031–2050. MR1329181

[34] POLITIS, D. N. and SHERMAN, M. (2001). Moment estimation for statistics from marked point processes. *J. R. Stat. Soc. Ser. B Stat. Methodol.* **63** 261–275. MR1841414

[35] POSSOLO, A. (1991). Subsampling a random field. In *Spatial Statistics and Imaging* (A. Possolo, ed.) 286–294. IMS, Hayward, CA. MR1195571

[36] SHERMAN, M. (1996). Variance estimation for statistics computed from spatial lattice data. *J. Roy. Statist. Soc. Ser. B* **58** 509–523. MR1394363

[37] SHERMAN, M. and CARLSTEIN, E. (1994). Nonparametric estimation of the moments of a general statistic computed from spatial data. *J. Amer. Statist. Assoc.* **89** 496–500. MR1294075

[38] SHORACK, G. R. (1982). Bootstrapping robust regression. *Comm. Statist. A—Theory Methods* **11** 961–972. MR0655465

[39] ZHU, J. and LAHIRI, S. N. (2006). Weak convergence of blockwise bootstrapped empirical processes for stationary random fields with statistical applications. *Statist. Inference Stochastic Processes.* To appear.
[40] ZHU, J., LAHIRI, S. N. and CRESSIE, N. (2002). Asymptotic inference for spatial CDFs over time. *Statist. Sinica* **12** 843–861. MR1929967
[41] ZHU, J. and MORGAN, G. D. (2004). Comparison of spatial variables over subregions using a block bootstrap. *J. Agric. Biol. Environ. Stat.* **9** 91–104.

DEPARTMENT OF STATISTICS  
IOWA STATE UNIVERSITY  
AMES, IOWA 50011  
USA  
E-MAIL: snlahiri@iastat.edu

DEPARTMENT OF STATISTICS  
UNIVERSITY OF WISCONSIN  
MADISON, WISCONSIN 53706  
USA  
E-MAIL: jzhu@stat.wisc.edu